\documentclass[12pt,draftcls,onecolumn]{IEEEtran}
\usepackage{amsfonts,amssymb}
\usepackage[cmex10]{amsmath}
\usepackage{epsfig}
\usepackage{amsfonts}
\usepackage{graphics}
\usepackage{graphicx}
\usepackage{subfig}
\usepackage{tabularx}
\usepackage{mathtools}
\usepackage{multirow}
\usepackage{cite}
\usepackage{psfrag}

\epsfverbosetrue

    \def\independenT#1#2{\mathrel{\setbox0\hbox{$#1#2$}%
    \copy0\kern-\wd0\mkern4mu\box0}}

\newcommand{\el}{{ ~~\it Example: \ }}

\title{Explicit Characterization of Stability Region for Stationary Multi-Queue Multi-Server Systems}  
\author{Hassan Halabian, Ioannis Lambadaris, Chung-Horng Lung
\\
Department of Systems and Computer Engineering\\
Carleton University, 1125 Colonel By Drive, Ottawa, ON, K1S 5B6, Canada\\
Email: \{hassanh, ioannis, chlung\}@sce.carleton.ca
}      
\begin{document}             
\maketitle                   
\newtheorem{theor}{Theorem}
\newtheorem{lem}{Lemma}
\newtheorem{pro}{Proposition}
\newtheorem{defin}{Definition}
\newtheorem{conj}{Conjecture}
\newtheorem{cor}{Corollary}

\IEEEpeerreviewmaketitle
\begin{abstract}
In this paper, we characterize the network stability region (capacity region) of multi-queue multi-server (MQMS) queueing systems with stationary channel distribution and stationary arrival processes. The stability region is specified by a finite set of linear inequalities. We first show that the stability region is a polytope characterized by the finite set of its facet defining hyperplanes. We explicitly determine the coefficients of the linear inequalities describing the facet defining hyperplanes of the stability region polytope.
We further derive the necessary and sufficient conditions for the stability of the system for general arrival processes with finite first and second moments. For the case of stationary arrival processes, the derived conditions characterize the system stability region. Furthermore, we obtain an upper bound for the average queueing delay of Maximum Weight (MW) server allocation policy which has been shown in the literature to be a throughput optimal policy for MQMS systems. Using a similar approach, we can characterize the stability region for a fluid model MQMS system. However, the stability region of the fluid model system is described by an infinite number of linear inequalities since in this case the stability region is a convex surface. We present an example where we show that in some cases depending on the channel distribution, the stability region can be characterized by a finite set of non-linear inequalities instead of an infinite number of linear inequalities.
\end{abstract}
\section{Introduction}
\IEEEPARstart{O}{ptimal} stochastic network control is one of the primary goals in the design of emerging wireless networks. One of the objectives of stochastic control in wireless networks is to enable cross layer designs to achieve stochastically optimal resource allocation in the physical and MAC layers, coupled with flow control/utility optimization strategies in transport layer and routing in the network layer.
 Examples of network resources at the MAC layer of wireless systems are OFDM subcarriers and CDMA codes and at the physical layer is transmission power.
Apart from the resource allocation problem in physical and MAC layers, flow control plays a crucial role in ensuring the system stability while achieving a level of network fairness. A flow control strategy must decide how much fraction of the injected traffic to the system must be admitted to assure the stability of queues in the network layer while achieving the optimal network fairness among users. A quantitative measure of fairness which is widely used in literature (e.g. \cite{phd,Now,berry2002, kelly1,kelly2,neely-infocom05,stolyar,boyd,atila}) is to define a set of utility functions $f_n(r)$ which illustrate the grade of satisfaction for each user $n$ while it transmits data traffic at rate $r$.
Consider a general network (wired or wireless) with $N$ source nodes. Suppose that $\lambda_n,~n=1,2,...,N$ be the traffic arrival rate of each user $n$. For such a network, the \textit{stability region} or \textit{network capacity region} is defined as the closure of the set of all arrival rate vectors for which there exists a resource allocation policy that can stabilize the system \cite{Now, phd, Tassiulas92}.
Let us denote the network stability region by $\Lambda$.
The flow control strategy determines the admitted rates $r_n$ from each user $n$ by solving the following flow maximization problem.
\begin{eqnarray}
\label{fairness}~~ \texttt{Maximize:}~~~~~~~~~~\sum_{n=1}^N f_n(r_n)~~~~~~~~~~~~~~~~~~~~  \\ 
\texttt{Subject to:}~~~~~~~~~~r=(r_1,r_2,...,r_N) \in \Lambda ~~~~~\nonumber \\
0 \leq r_n \leq \lambda_n~~\forall n=1,2,...,N \nonumber\hspace{-.5cm}
\end{eqnarray}
By solving the above maximization problem, each node adjusts its admitted traffic rate to the computed optimal point from (\ref{fairness}). We assume that functions $f_n(r_n)$ are non-decreasing and concave. The choice of functions $f_n(r_n)$ depends on the desired fairness properties in the network. For example, choosing the utility function $f(r)=\log (1+\beta r)$ (for some large constant parameter $\beta>0$) will result in proportional fairness behaviour \cite{kelly1,Now}.

To solve the optimization problem in (\ref{fairness}), we need to have the \textit{stability region} of the system.
Note that the stability region is unique for each network and independent of the resource allocation policy. Although stability region has been described and characterized for a general network in \cite{Tassiulas92,Tassiulas97,Now,phd} as given by the convex hull of a set of fixed points, such an \textit{implicit characterization} cannot be directly used to define the constraints of network optimization problem (\ref{fairness}). The main goal of this paper is to introduce a linear algebraic characterization of the stability region for stationary Multi-Queue Multi-Server (MQMS) queueing systems. Such queueing systems can be used to model practical multi-user wireless networks with multiple orthogonal sub-channels such as OFDM sub-carriers. In such networks, there is a set of users generating random packet arrivals and a set of shared orthogonal sub-channels (servers) that are assigned to 
users according to certain rules.
Because of users mobility, environmental changes, fading, etc., the channel quality of each user to each server is changing randomly with time.
Therefore, resource allocation in such networks can be modeled as a server allocation problem in multi-queue multi-server queueing systems with time varying channel conditions \cite{javidi,javidi2,javidi3,javidi4, hussein, ganti,khodam}. Our focus in this paper would be on the server allocation problem in MQMS systems with stationary channel distribution for which we will introduce a linear algebraic representation of the network stability region (i.e., $\Lambda$). Specifically, we will determine explicitly all the coefficients of the linear inequalities that describe the stability region. These inequalities then can be tabulated and used as the constraints of flow optimization problems similar to (\ref{fairness}).

The stability problem in wireless queueing networks was mainly addressed in \cite{Tassiulas92, Tassiulas93,Tassiulas97,Now, phd}. In \cite{Tassiulas92}, authors introduced the notion of stability region of a queueing network. They considered a time slotted system in their work and assumed that arrival processes are i.i.d. sequences and the queue length process is a Markov process. In \cite{Tassiulas93}, they characterized the network stability region of multi-queue single-server systems with time varying ON-OFF connectivities. They also proved that for a symmetric system (with the same arrival and connectivity statistics for all the queues), LCQ (Longest Connected Queue) policy maximizes the stability region and also provides the optimal performance in terms of average queue occupancy (or equivalently average queueing delay). In \cite{Now,phd} and \cite{power}, the notion of network stability region of a wireless network was introduced for more general arrival and queue length processes. Furthermore, Lyapunov drift techniques were applied in \cite{phd} and \cite{Now} to analyze the stability of the proposed policies for stochastic optimization problems in wireless networks.

The problem of server allocation in multi-queue multi-server systems with time varying connectivities was mainly addressed in \cite{javidi,javidi2,javidi4, hussein, khodam}. In \cite{javidi2}, Maximum Weight (MW) policy was proposed as a throughput optimal server allocation policy for MQMS systems with stationary channel process. However, in \cite{javidi2} the conditions on the arrival traffic to guarantee the stability of MW were not explicitly mentioned. 
In our previous work in \cite{khodam}, we characterized the network stability region of multi-queue multi-server systems with time varying ON-OFF channels. We also obtained an upper bound for the average queueing delay of AS/LCQ (Any Server/Longest Connected Queue) policy which is the throughput optimal server allocation policy for such systems.

References \cite{javidi,javidi3,javidi4,hussein} study the optimal server allocation problem in terms of average queueing delay. In \cite{javidi,javidi3,javidi4}, the authors argue that in general, achieving instantaneous throughput and load balancing is impossible in a general MQMS system. However, they showed that this goal is attainable in the special case with ON-OFF channel processes. They also introduced MTLB (Maximum-Throughput Load-Balancing) policy and showed that this policy is minimizing a class of cost functions including total average delay for the case of two symmetric queues. The work in \cite{hussein} considers this problem for general number of symmetric queues and servers. Authors in \cite{hussein} characterized a class of \textit{Most Balancing} (MB) policies among all work conserving policies which are minimizing, in stochastic ordering sense, a class of cost functions including total average delay. They used stochastic ordering and dynamic coupling arguments to show the optimality of MB policies for symmetric systems.

In this paper, we will characterize the stability region of multi-queue multi-server queueing system with stationary channel and arrival processes. Toward this, the necessary and sufficient conditions for the stability of the system are derived for general arrival processes with finite first and second moments. For stationary arrival processes, these conditions establish the network stability region of such systems. 
Our contribution in this work is to characterize the stability region as a \textit{convex polytope} specified by a \textit{finite set of linear inequalities} that can be numerically tabulated and used to solve network optimization problems similar to (\ref{fairness}) (Refer to Lemmas \ref{l1} to \ref{l5} and Theorems \ref{t1} and \ref{sc}).
%
Later in the paper, we further introduce an upper bound for the average queueing delay of MW policy \cite{javidi2} which is a throughput optimal policy for MQMS systems. We also study the stability of fluid model MQMS systems for which using the same approach as what we use in the packetized system, we characterize the linear algebraic representation of the stability region. We show that the stability region in this case is characterized by an infinite number of linear inequalities. However, by an example we show that depending on the channel distribution and the dimension of the system, we may characterize the stability region by a limited number of non-linear inequalities instead of infinite number of linear inequalities.  

The rest of this paper is organized as follows. In section \ref{notation}, we introduce the notation required throughout this paper. Section \ref{modeldes} describes the queueing model we focused on. In section \ref{back}, we discuss about the notion and definition of \textit{strong stability} in queueing networks. We also briefly elaborate on the Lyapunov drift technique used to prove the system stability \cite{phd,Now}. Moreover, we review some fundamental properties of polytopes in section \ref{back}. In section \ref{sec-asli}, we will derive necessary and sufficient conditions for the stability of our model. We determine the coefficients of the linear inequalities describing the facets of the stability region polytope. We also find an upper bound for the average queue occupancy (or average queueing delay). Finally, we introduce fluid model MQMS systems and study the stability region of such queueing systems. Section \ref{conc} presents the conclusions of our work. 
\section{Notation} \label{notation}
In this section, we introduce basic notation used throughout the paper. Additional notation will be introduced when necessary. All the vectors are considered to be row vectors. By $\underline{1}_K$ ($\underline{0}_K$), we denote a row vector of size $K$ whose elements are all identically equal to $``1"$ ($``0"$).
The time average of a function $f(t)$ is denoted by $ \overline{f(t)} $, i.e.,  $\overline{f(t)}= \lim_{t\rightarrow \infty}\frac{1}{t} \sum_{\tau=1}^{t} f(\tau)$. The operator $``\circledast"$ is used for entry-wise multiplication of two matrices. 
The expectation of random processes (or random variables) is denoted by $E[\cdot]$. The cardinality of a set is denoted by $|\cdot|$. The operator for inner product of two vectors is $\left\langle \cdot, \cdot \right\rangle$. The boundary of a set is represented by $bound (\cdot)$. For any vector $\alpha=(\alpha_1,\alpha_2,...,\alpha_N)$ and a non-empty ordered (sub)set of indices $ \mathcal{U}=\{u_1,u_2,...,u_{|\mathcal{U}|}\} \subseteq \{1,2,...,N\} $ and $u_1<u_2<...< u_{|\mathcal{U}|}$, we define $ \alpha_{\mathcal{U}}=(\alpha_{u_1},\alpha_{u_2},...,\alpha_{u_{|\mathcal{U}|}}) $.

\section {Model Description} \label{modeldes}
\begin{figure}[tp]
    \centering
    \psfrag{a}[][][.9]{$C_{1,1}(t)$}
    \psfrag{b}[][][.9]{$C_{1,2}(t)$}
    \psfrag{c}[][][.9]{$C_{1,K}(t)$}
    \psfrag{d}[][][.9]{$C_{N,1}(t)$}
    \psfrag{e}[][][.9]{$C_{N,2}(t)$}
    \psfrag{f}[][][.9]{$C_{N,K}(t)$}  
    \psfrag{g}[][][1]{$X_1(t)$}
    \psfrag{h}[][][1]{$X_2(t)$}
    \psfrag{i}[][][1]{$X_N(t)$}
    \psfrag{j}[][][1]{$A_1(t)$}
    \psfrag{k}[][][1]{$A_2(t)$}
    \psfrag{l}[][][1]{$A_N(t)$}
    \includegraphics[width=0.45\textwidth]{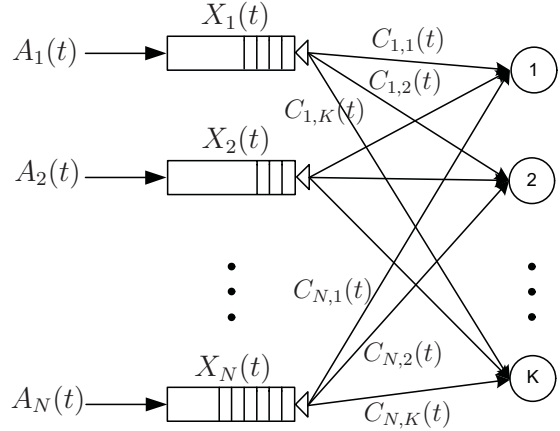}
	\caption{ Multi-queue multi-server queueing system with stationary channel distribution}
	\label{model} \vspace*{-5mm}
\end{figure} 
We consider a time slotted queueing system with equal length time slots and equal length packets. The model consists of a set of parallel queues ${\cal N}=\{1,2,...,N\}$ and a set of identical servers ${\cal K}=\{1,2,...,K\}$. Each server can serve at most one queue at each time slot, i.e., we do not allow server sharing in the system. At each time slot $t$, the capacity of the link between each queue $ n\in \mathcal{N}$ and server $k\in \mathcal{K}$ is assumed to be $C_{n,k}(t)$ packets/time slot (Figure \ref{model}), i.e., during time slot $t$, server $k$ can serve at most $C_{n,k}(t)$ packets of queue $n$ successfully if it is allocated to queue $n$ at that time slot. We assume that $C_{n,k}(t) \in \mathcal{M}$ where $\mathcal{M}=\{m\in \mathbb{Z}_+ \mid m \leq M\} $, for a given $M$. Therefore, at each time slot $t$, the channel state may be expressed by an $N\times K$ matrix $C(t)=\left( C_{n,k}(t)\right),n\in \mathcal{N}, k \in \mathcal{K}, C_{n,k}(t)\in \mathcal{M}$. The channel process is defined as ${\{C(t)\}}_{t=1}^{\infty}$ with the state space $\cal S$.
Note that $\cal S$ is a finite set with $|\mathcal{S}|= (M+1)^{NK}$. We will label each element of $\cal S$ by a positive integer index $s \in \{1,2, ... ,(M+1)^{NK}\}$. Suppose that the channel state matrix associated to the channel state $s$ is denoted by $C_{s}$. The channel process is assumed to have stationary distribution with stationary probabilities 
$\pi_s=\Pr(C(t)=C_s)$.

Queues are fed by exogenous packet arrival processes $A_n(t),n=1,2,...,N$,
i.e., the number of packet arrivals to queue $n$ during time slot $t$ is represented by $A_n(t)$. The arrival vector at time slot $t$ is denoted by $A(t)=(A_1(t),A_2(t),...,A_N(t))$. For these processes, suppose that $E[A_n^2(t)]\leq A_{max}^2<\infty$ \textit{for all} $t$. We assume that each queue has an infinite buffer space. We also assume that new arrivals are added to each queue at the end of each time slot. Let $X(t)=(X_1(t), ... , X_N(t))$ be the queue length vector at the end of time slot $t$ after adding new arrivals to the queues\footnote{We assume that at each time slot departures occur first and then at the end of the time slot new arrivals are added to the queues.}.

A server scheduling policy at each time slot should decide how to allocate servers from set $\cal K$ to the queues in set \(\cal N\). This must be accomplished based on the available information about the channel state of the system at time slot $t$ (i.e., $C(t)$) and also the queue length state at the beginning of time slot $t$ (i.e., $X(t-1)$). Therefore, at the beginning of each time slot $t$ the scheduler has to determine an allocation matrix $I(t) \in \mathcal{I}$ where
\[\mathcal{I}=\left\{I_{N\times K}=\left(I_{n,k}\right),I_{n,k}\in \{0,1\} \mid ~ \sum_{n=1}^N I_{n,k}\leq 1 ~\forall k\in \mathcal{K} \right\} .\]
and $\cal I$ is called allocation matrix space. We observe that each matrix in $\cal I$ can have at most a single $``1"$ in each of its columns.
The queue length vector evolves with time according to the following rule.
\begin{equation}
X^{\mathrm{T}}(t)=\left(X^{\mathrm{T}}(t-1)-\left( C(t)\circledast I(t) \right)\underline{1}_K^{\mathrm{T}} \right)^+ +A^{\mathrm{T}}(t)
\end{equation}
where $(\cdot)^+$ is defined as follows: For an arbitrary vector $v$ of size $|v|$, $(v)^+$ is a vector of the same size whose $i$'th element, $(v)^+_i$, is defined by
\begin{equation}
(v)^+_i= \left\{
  \begin{array}{l l}
    0 & \quad \text{if $v_i<0$ }\\
     v_i & \quad \text{if $v_i\geq 0$}\\
  \end{array} \right. .
\end{equation}

\section{Background} \label{back}

\subsection{Strong Stability Definition}     
We now introduce the definition of strong stability for a queueing system \cite{Now, phd}. Other definitions can be found in \cite{Asmussen, Tassiulas92, Tassiulas93, Leonardi}.
Consider a discrete time single queue system with an arrival process \(A(t)\) and service process $\mu (t)$. As we mentioned earlier, the arrivals are added to the system at the end of each time slot. We can see that the queue length process $X(t)$ evolves with time according to the following recursion;
\begin{equation} \label{sq}
X(t)=(X(t-1)-\mu(t))^+ +A(t)
\end{equation}
\begin{defin} 
A queue evolving with time according to (\ref{sq}) is said to be \emph{strongly stable} \cite{Now} if 
\begin{equation}
\limsup_{t \to \infty} {1\over{t}} \displaystyle\sum _ {\tau =0}^{t-1}E[X(\tau)] < \infty.
\end{equation}   
\end{defin} 

\begin{defin}  
A queueing network is said to be strongly stable \cite{Now} if all the queues in the system are strongly stable. 
\end{defin} 

In this paper, we will employ the definition of 
strong stability and in what follows we use the terms ``stability" and ``strong stability" interchangeably. In \cite{Now,phd}, it was proved that if a queue is strongly stable and if for all $t$ either $E[A(t)] \leq A$ or $E[\mu (t)-A(t)] \leq D$ where $A$ and $D$ are finite given non-negative constants, then 
\begin{equation} \label{zerolim}
\lim _{t \to \infty} {1\over{t}} E[X(t)] =0.
\end{equation}

\subsection{Lyapunov Drift}
A very important and useful mathematical tool used in network stability analysis and stochastic control/optimization of wireless networks is the \emph{Lyapunov Drift} technique introduced in \cite{Now,phd,Leonardi, McKeown,power}. The main idea behind the Lyapunov stability method is to define a nonnegative function of queue backlogs which can be seen as a \emph{measure} of the total aggregated backlog in the system at time slot \(t\). Then, we evaluate the ``drift" of such function in two successive time slots by taking the effect of control decisions (scheduling or resource allocation policy) into account. If the expected value of the drift is negative as the backlog goes beyond a fixed threshold, then the system is stable. This method was used in \cite{Now,phd,McKeown,power,  Leonardi} to prove the stability of a number of queueing systems.

For a queueing system with $N$ queues and queue length vector \(X(t)=(X_1(t), ... , X_N(t))\), the following quadratic function has been usually used in literature (\cite{Now, phd, McKeown, power,Leonardi}) as a Lyapunov function.
\begin{equation}
\label{lypfunc}L(X)= \displaystyle\sum_{n=1}^N X_n^2(t)
\end{equation}

Assume that $E[X_n(0)]<\infty$, $\forall n \in \mathcal{N}$ and $X(t)$ evolves with some probabilistic law (not necessarily Markovian). Then, the following holds \cite{Now}.
\begin{lem}
\label{llyap}
If there exist constants \(B>0\) and \(\delta >0\) such that for all time slots \(t\) we have
\begin{eqnarray} 
 \label{lyapunov} E[L(X(t+1))-L(X(t)) \mid X(t)] \leq B - \delta\displaystyle\sum _{n=1}^N X_n(t),
\end{eqnarray}
then the system is strongly stable and further we have
\begin{equation}
\limsup_{t \to \infty} {1\over{t}} \displaystyle\sum_ {\tau =0}^{t-1} \displaystyle\sum _{n=1}^N E[X_n(\tau)] \leq \frac{B}{\delta}.
\end{equation}
\end{lem}

The proof of the lemma can be found in \cite{phd,Now}. The left hand side of expression (\ref{lyapunov}) is usually called Lyapunov drift function which is a measure of the expected value of changes in the backlog in two successive time slots. From Lemma \ref{llyap}, we can easily see the idea behind Lyapunov method in stabilizing queueing systems. It is not hard to show that, when the aggregated backlog in the system goes beyond the bound \({B} \over \delta\), the Lyapunov drift in the left hand side of (\ref{lyapunov}) will be negative, meaning that the system receives a negative drift on the expected aggregated backlog in two successive time slots. In other words, the system tends toward lower backlogs and this results in its stability.

\subsection{Fundamental Concepts of Polytopes} \label{polytopesconcepts}
We present a brief review on convex polytopes and fundamental properties of them.
These concepts will be needed for specifying the stability region of MQMS system in section \ref{sec-asli}.



\begin{defin} 
\label{polyt}
A convex polytope is defined as the convex hull of a finite set of points \cite{lee,alex}.
\end{defin}
%

According to \textit{Weyl's Theorem} \cite{lee}, a polytope in $\mathbb{R}^N$ always can be expressed by a set $\mathcal{P}=\left\{x \in \mathbb{R}^N  \mid \alpha_{\ell} x^{\mathrm{T}}\leq \beta_{\ell}~for~\ell=1,2,...,L  \right\}$ for some positive integer $L$ and $\alpha_{\ell} \in \mathbb{R}^N$ and $\beta_{\ell} \in \mathbb{R}$.

\begin{defin}
The dimension of a polytope $\cal P$ is represented by $\dim (\mathcal{P})$ and is defined as one less than the maximum number of affinely independent points in $\cal P$ \cite{lee,alex}.
\end{defin}
A polytope $\mathcal{P} \subset \mathbb{R}^N$ is said to be full dimensional if $\dim (\mathcal{P})=N$. 

\textit{Dimension Theorem}: For a polytope $\mathcal{P} \subset \mathbb{R}^N$, dimension of $\cal P$ is equal to $N$ minus the maximum number of linearly independent equations satisfied by all the points in $\cal P$. 
\begin{defin}
For given $a$ and $b$, equality $ax^{\mathrm{T}}\leq b$ is called \textit{valid} for polytope $\mathcal P$ if for every point $x_0 \in \mathcal P$, $ax_0^{\mathrm{T}}\leq b$.
\end{defin}
\begin{defin}
A \textit{face} of polytope $\mathcal P$ is defined as $\mathcal{F}=\{x \in \mathcal{P} \mid ax^{\mathrm{T}}= b\}$ where 
inequality $ax^{\mathrm{T}}\leq b$ is a valid inequality for $\mathcal P$.
\end{defin}

We call the valid inequality $ax_0^{\mathrm{T}}= b$ a \textit{face defining hyperplane} for $\mathcal{P}$ if its associated face is not empty. Therefore, $ax_0^{\mathrm{T}}= b$ is a face defining hyperplane for $\mathcal{P}$ if it intersects with $\mathcal{P}$ at least at one point.
Note that $\mathcal{P}$ has finitely many faces. However, the face defining hyperplanes of a polytope can be infinite.
\begin{defin} 
\label{facet}
A \textit{facet} of polytope $\mathcal P$ is a maximal face distinct\footnote{Maximal relative to inclusion} from $\cal P$ \cite{alex}. Hence, all faces of $\mathcal{P}$ with dimension $\dim (\mathcal{P})-1$ are called facets of $\mathcal{P}$.
\end{defin}

For a polytope $\mathcal{P}=\left\{x \in \mathbb{R}^N  \mid \alpha_{\ell} x^{\mathrm{T}}\leq \beta_{\ell}~for~\ell=1,2,...,L  \right\}$ an inequality is \textit{redundant} if polytope $\cal P$ remains unchanged by removing the inequality.

\textit{Redundancy Theorem in Polytopes} \cite{lee}: Face defining hyperplanes describing faces of dimension less than $\dim (\mathcal{P})-1$ are redundant.

Redundancy theorem states that to describe a polytope completely, only facet defining hyperplanes are sufficient.

\section{Linear Algebraic Representation of Stability Region Ploytope} \label{sec-asli}
In this section, we first introduce the convex hull representation of the stability region \cite{Now,phd}. Then, we derive the necessary and sufficient conditions for the stability of MQMS system with stationary channel processes and general arrival processes with finite first and second moments. It will be shown that for stationary packet arrival processes, these conditions establish the stability region polytope. We determine a linear algebraic representation of the stability region in which we explicitly determine the coefficients of the linear inequalities describing the stability region.
\vspace{-2mm}
\subsection{Stability Region Geometry}
Consider the class of deterministic policies $\cal G$ where at each state of the system, each policy $g\in \cal G$ allocates the servers according to a predetermined allocation matrix depending on the channel state matrix (i.e., matrix $C(t)$). More specifically, for each policy $g \in \cal G$, there exists a one to one mapping from the channel state space $\cal S$ to the allocation matrix space $\cal I$, namely $I^{(g)}: \mathcal{S} \longmapsto \mathcal{I}$.
Therefore, each deterministic policy $g$ is specified by $|\cal S|$ allocation matrices $I_s^{(g)}, ~s\in \cal S$. The scheduler observes the state of the system and then based on the observed state $s$, it allocates the servers according to an allocation matrix $I^{(g)}_s$. Note that set $\cal I$ is a finite set and since the channel state space is also finite, set $\cal G$ is finite.

Each deterministic policy $g$ provides an average transmission rate for each queue $n$. Let $R_n^{(g)}$ denote the time averaged transmission rate provided to queue $n$ and $R^{(g)}=(R_1^{(g)},R_2^{(g)},...,R_N^{(g)})$ the vector of average transmission rates. By conditioning on the channel state of the system, it is not hard to see that for each deterministic server allocation policy $g$, we have
\begin{eqnarray}
\label{rg}
R^{(g)}={\left(\displaystyle\sum_{s\in {\cal S}}\pi_s \left(C_s\circledast I_s^{(g)} \right)\underline{1}_K^{\mathrm{T}}\right)}^{\mathrm{T}}.
\end{eqnarray}   

Note that each rate vector $R^{(g)}$ determines a single point in $\mathbb{R}_{+}^N$. Now, consider the convex hull of all the points $R^{(g)}$, $g \in \cal G$ in $\mathbb{R}_{+}^N$, i.e.,   
\begin{eqnarray}
\label{conv}
\mathcal{R}=conv.hull_{\substack{~\\~\\ \hspace{-16 mm} g \in \mathcal{G}}}(R^{(g)}).
\end{eqnarray}

Each point in $\cal R$, let say $R^{\star}$ can be represented by a convex combination of a finite set of points, $R^{(g)},g\in \cal G$, i.e., 
\begin{eqnarray}
\label{conv.comb}
R^{\star}=\displaystyle\sum_{i=1}^{|\mathcal{G}|}p_i R^{g_i}~,g_i\in \mathcal{G}~,~\displaystyle\sum_{i=1}^{|\mathcal{G}|}p_i=1,~p_i\geq 0.
\end{eqnarray}
Hence, $\cal R$ contains all the achievable transmission rate vectors of MQMS system.
To achieve the transmission rate vector $R^{\star} \in \mathcal{R}$ it is enough to select policy $g_i$ with probability $p_i$, i.e., in $p_i$ fraction of time slots. In other words, all the transmission rate vectors $R^{\star} \in \cal R$ are achievable by applying a \textit{randomized} policy that at each time slot selects policy $g_i$ with probability $p_i$. Thus, the following Lemma follows.

\begin{lem}
\label{l2}
The set of achievable transmission rate vectors $\cal R$ is specified by a polytope.
\end{lem}

\begin{IEEEproof}
 The lemma follows directly from the definition of polytope in (\ref{polyt}), equation (\ref{conv}) and the fact that set $\cal G$ is finite. 
\end{IEEEproof}
We will denote the achievable transmission rate polytope by $\cal P$. According to the definition of stability region and the discussion in \cite{Now, phd,Tassiulas97, power} regarding the network stability region, we can conclude that polytope $\cal P$ specifies the stability region of MQMS queueing system with stationary channel distribution and stationary arrival processes (in this case $E[A_n(t)]=\lambda_n$ and $\lambda=(\lambda_1,\lambda_2,...,\lambda_N)$). Specifically, it has been shown that if the system is stable, we should have $\lambda \in \cal P$ and also if $\lambda \in {\cal P} -bound({\cal P})$, then there exists a server allocation policy that can stabilize the system. 

Although stability region of MQMS system is described by (\ref{conv.comb}), it cannot be applied as the constraints of utility optimization problem (\ref{fairness}). In fact, (\ref{conv.comb}) only provides us an implicit description of stability region. To solve network optimization problems like (\ref{fairness}), we should characterize the stability region by a set of linear/non-linear convex inequalities (or equalities) which is given in the following. 

\vspace{-2mm}

\subsection{Necessary Conditions for the Stability of MQMS System}
We first give an outline of the steps we take to find the necessary conditions for the stability of the system.
\begin{itemize}
\item  Using Lemmas \ref{l1}-\ref{l4}, we will characterize all the face defining hyperplanes of the stability region polytope. More specifically, we will determine the normal vector associated to each hyperplane as well as a single point where the hyperplane touches the stability region polytope. We show that each vector in $\mathbb{R}_+^N$ is associated to a face defining hyperplane.
\item The stability region is a polytope and therefore can be characterized by the finite set of its facet defining hyperplanes. We will introduce a subset of $\mathbb{R}_+^N$, namely $V\subset \mathbb{R}_+^N$ and in Lemma \ref{l6} we show that the vectors outside set $V$ cannot be associated to a facet defining hyperplanes of the stability region polytope. Note that set $V$ is an infinite set itself.
\item In Lemma \ref{l5}, we will show that although set $V$ is an infinite set, it is finite up to multiplication of vectors by positive scalars, i.e., all the vectors in $V$ can be produced by multiplying a positive scalar to a vector in a finite set $\widehat{V}$.
\item Using Lemmas \ref{l1} to \ref{l5}, we will prove Theorem \ref{t1} which states the necessary conditions for the stability of the system. We also argue that the derivation of set $\widehat{V}$ becomes difficult for large $N$ and $M$. We introduce a finite superset of set $\widehat{V}$ whose elements can be specified easily. 
\end{itemize}

We introduce the departure matrix $H_{N\times K}(t)=(H_{n,k}(t)), n\in {\cal N}, k\in{\cal K}$ in which $H_{n,k}(t)$ ($H_{n,k}(t)\leq C_{n,k}(t)$) represents the total number of packets served by server $k$ from queue $n$ at time slot $t$. 
Thus, the departure process from queue $n$ at time slot $t$ would be $\sum_{k=1}^K H_{n,k}(t)$.
The following equality illustrates the evolution of queue length process with time.
\begin{eqnarray} 
\label{evol} X_n(t)=X_n(t-1)-\displaystyle\sum _{k=1}^K H_{n,k}(t) +A_n(t)~~~~~~\forall n \in \cal N
\end{eqnarray}

For a strongly stable MQMS queueing system we can prove the following Lemma. 

\begin{lem}
\label{l1}If the MQMS system is strongly stable, then for any vector $\alpha=(\alpha_1,\alpha_2,...,\alpha_N)\in \mathbb{R}^N$ we have
\begin{eqnarray}
\label{l1_1} \overline{\alpha E[A(t)]^{\mathrm{T}}}=\overline{\alpha E[H(t)]\underline{1}_K^{\mathrm{T}}}.
\end{eqnarray}

\end{lem}

\begin{IEEEproof}
The recursion (\ref{evol}) will result into
\begin{eqnarray}
\label{l1_2}
X_n(t)=X_n(0)-\displaystyle\sum _{\tau=1}^t \displaystyle\sum_{k=1}^K H_{n,k}(\tau) + \displaystyle\sum _{\tau=1}^t A_n(\tau).
\end{eqnarray}
By multiplying vector $\alpha$ to the queue length vector $X(t)$ we have
\begin{eqnarray}
\label{l1_3}
\alpha X^{\mathrm{T}}(t)=\alpha X^{\mathrm{T}}(0) - \displaystyle\sum _{\tau=1}^t \alpha H(\tau)\underline{1}_K^{\mathrm{T}} + \displaystyle\sum _{\tau=1}^t \alpha A^{\mathrm{T}}(\tau).
\end{eqnarray}
Taking the expectation from both sides, dividing by \(t\) and then taking the limit as \(t\) goes to infinity, we will have the following.
\begin{eqnarray}
\label{l1_4}
\lim_{t\rightarrow \infty}\frac{1}{t}\alpha E[ X(t)]^{\mathrm{T}}=\lim_{t\rightarrow \infty}\frac{1}{t}\alpha E[X(0)]^{\mathrm{T}} -  \lim_{t\rightarrow \infty} \frac{1}{t} \displaystyle\sum_{\tau=1}^t \alpha E[H(\tau)]\underline{1}_K^{\mathrm{T}} + \lim_{t\rightarrow \infty} \frac{1}{t} \displaystyle\sum _{\tau=1}^t \alpha E[A(\tau)]^{\mathrm{T}}
\end{eqnarray}
According to (\ref{zerolim}) and the assumption that \(E[X(0)]<\infty\), the left hand side term and the first term in the right hand side term are equal to zero and therefore the result follows.
\end{IEEEproof}

A direct result of Lemma \ref{l1} for a single queue $n$ is the following.
\begin{eqnarray} 
\label{a=h}
\overline{E[A_n(t)]}=\overline{E\left[\sum_{k=1}^K H_{n,k}(t)\right]}
\end{eqnarray}
In other words, for a strongly stable queue $n$, the time averaged total expected arrival to queue $n$ is equal to the time averaged total expected departure from queue $n$.


\begin{lem}
\label{l3}
If the MQMS system is strongly stable, then  $\overline{E[A(t)]} \in \cal P$.
\end{lem}

\begin{IEEEproof}
The proof follows directly by considering (\ref{a=h}) and the fact that 
\begin{eqnarray}
\label{h<r} \forall~ H(t),~\exists R \in \mathcal{P} : \overline{E\left[\sum_{k=1}^K H_{n,k}(t)\right]}\leq R_n~~ \forall n\in \cal N .
\end{eqnarray} 
\end{IEEEproof}

Based on Lemmas \ref{l1} and \ref{l3}, we can prove the following Lemma.
\begin{lem}
\label{l4}
If the MQMS system is strongly stable, then \textit{for all} $\alpha \in \mathbb{R}_+^N$
\begin{eqnarray} 
\label{l4_1} 
\alpha \overline{E[A(t)]^{\mathrm{T}}} \leq \displaystyle\sum_{s\in \mathcal{S}}\pi_s 
  \max_{I \in \mathcal I}  \left(\alpha (C_s \circledast I )\underline{1}_K^{\mathrm{T}} \right)
\end{eqnarray}
\end{lem}

\begin{IEEEproof}
Consider $\alpha \in \mathbb{R}^N$. Since the system is strongly stable, from Lemma \ref{l1} we have
\begin{eqnarray}  
\label{l4_2}
\lim_{t\rightarrow \infty}\frac{1}{t} \displaystyle\sum_{\tau=1}^{t} \alpha E[A(\tau)]^{\mathrm{T}} =  \lim_{t\rightarrow \infty} \frac{1}{t}\displaystyle\sum_{\tau=1}^{t} \alpha E[H(\tau)]\underline{1}_K^{\mathrm{T}} 
\end{eqnarray}
By conditioning on the channel state process, we will have
\begin{eqnarray}  
\label{l4_3}
\lefteqn{\lim_{t\rightarrow \infty} \frac{1}{t}\displaystyle\sum_{\tau=1}^{t} \alpha E[H(\tau)]\underline{1}_K^{\mathrm{T}} }\nonumber\\
&& = \lim_{t\rightarrow \infty} \frac{1}{t} \displaystyle\sum_{\tau=1}^{t} \displaystyle\sum_{s\in \mathcal{S}}\pi_s \alpha E[H(\tau) \mid S(\tau)=s]\underline{1}_K^{\mathrm{T}} \nonumber \\
&& \leq  \lim_{t\rightarrow \infty} \frac{1}{t} \displaystyle\sum_{\tau=1}^{t} \displaystyle\sum_{s\in \mathcal{S}}\pi_s 
 (\alpha)^+ E\left[  C_s\circledast I(\tau) \right]\underline{1}_K^{\mathrm{T}}  \nonumber \\
&& \leq  \lim_{t\rightarrow \infty} \frac{1}{t} \displaystyle\sum_{\tau=1}^{t} \displaystyle\sum_{s\in \mathcal{S}}\pi_s 
  \max_{I \in \mathcal I}  \left((\alpha)^+ (C_s \circledast I )\underline{1}_K^{\mathrm{T}} \right)  \nonumber \\
&& = \displaystyle\sum_{s\in \mathcal{S}}\pi_s 
  \max_{I \in \mathcal I}  \left( (\alpha)^+ (C_s \circledast I )\underline{1}_K^{\mathrm{T}} \right)
\end{eqnarray}
Now, consider two different cases:
\begin{itemize}
\item $\alpha \in \mathbb{R}^N_+$: The result follows directly from (\ref{l4_2}) and (\ref{l4_3}) since $(\alpha)^+=\alpha$. 
\item $\alpha \notin \mathbb{R}^N_+$: In this case, from (\ref{l4_2}) and (\ref{l4_3}) we have the following inequality. 
\begin{eqnarray} \label{j1}
\alpha \overline{E[A(t)]^{\mathrm{T}}} \leq \sum_{s\in \mathcal{S}}\pi_s 
  \max_{I \in \mathcal I}  \left((\alpha)^+ (C_s \circledast I )\underline{1}_K^{\mathrm{T}} \right)
\end{eqnarray}
However, since $(\alpha)^+ \in \mathbb{R}^N_+$ and according to the previous case we also have \begin{eqnarray} \label{j2}
(\alpha)^+ \overline{E[A(t)]^{\mathrm{T}}} \leq \sum_{s\in \mathcal{S}}\pi_s 
  \max_{I \in \mathcal I}  \left((\alpha)^+ (C_s \circledast I )\underline{1}_K^{\mathrm{T}} \right).
\end{eqnarray}
Noting (\ref{j2}) and the fact that $E[A(t)] \in \mathbb{R}^N_+$, we conclude that (\ref{j1}) is a redundant inequality. 
\end{itemize} 
\end{IEEEproof}


Since $\alpha \in \mathbb{R}_+^N$, the set of inequalities in (\ref{l4_1}) forms an infinite set. Each inequality in (\ref{l4_1}) determines a valid inequality for polytope $\cal P$.
 
\textbf{Fact}: The hyperplanes associated to the valid inequalities of (\ref{l4_1}) are all face defining hyperplanes of polytope $\cal P$.

To show this fact, let ${\cal I}^{\alpha}$ denote a set of allocation matrices $\{I^{\alpha}_s, s\in \cal S \}$ that maximize the right hand side of (\ref{l4_1}), i.e.,
\begin{eqnarray}
\label{argmax}I_s^{\alpha} = \arg\max_{I \in \cal I} \alpha (C_s \circledast I )\underline{1}_K^{\mathrm{T}}.
\end{eqnarray}
${\cal I}^{\alpha}$ is not unique and there may be more than one set of allocation matrices of ${\cal I}^{\alpha}$ whose elements maximize  $\alpha (C_s \circledast I )\underline{1}_K^{\mathrm{T}}$.
According to (\ref{rg}) and (\ref{conv}), ${\left(\sum_{s\in \mathcal{S}}\pi_s 
\left( C_s \circledast I_s^{\alpha} \right) \underline{1}_K^{\mathrm{T}}\right)}^{\mathrm{T}} \in \cal P$. On the other hand, since 
$\sum_{s\in \mathcal{S}}\pi_s 
 \alpha (C_s \circledast I^{\alpha}_s)\underline{1}_K^{\mathrm{T}} =\alpha \sum_{s\in \mathcal{S}}\pi_s 
  (C_s \circledast I^{\alpha}_s )\underline{1}_K^{\mathrm{T}},
$
 point  $ \sum_{s\in \mathcal{S}}\pi_s 
  (C_s \circledast I^{\alpha}_s )\underline{1}_K^{\mathrm{T}} $ is located on the hyperplane associated to (\ref{l4_1}).
Therefore, the set of inequalities in (\ref{l4_1}) determines all the non-empty faces of polytope $\cal P$.
To clarify what we proved in Lemmas \ref{l1} to \ref{l4} consider the following example.

\begin{figure}[h]
    \centering
    \psfrag{a}[][][1]{$\lambda_1$}
    \psfrag{b}[][][1]{$\lambda_2$}
    \psfrag{c}[][][1]{$X_1(t)$}
    \psfrag{d}[][][1]{$X_2(t)$}
    \psfrag{s}[][][1]{$S$}
    \psfrag{e}[][][1]{$p_1$}  
    \psfrag{f}[][][1]{$p_2$}
    \includegraphics[width=0.35\textwidth]{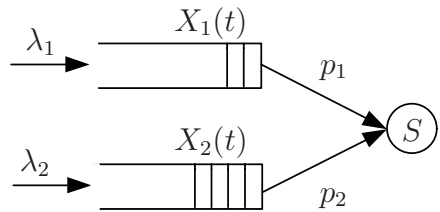}
	\caption{Two-queue single-server example} \vspace{-3mm}
	\label{example}
\end{figure}

\el Consider a system with $N=2$, $M=1$ (ON-OFF channels) and $K=1$ (see Figure \ref{example}). Assume that queues $1$ and $2$ are connected to the server with probabilities $p_1$ and $p_2$, respectively. This system with $N\geq 2$ queues was studied in \cite{Tassiulas93} when studying dynamic server allocation to a set of parallel queues. Here we consider a system with just $2$ queues so that we can illustrate the stability region in two dimensions. For such a system it was proven that the stability region is described by the following set of inequalities. 
\begin{gather}
\lambda_1\leq p_1\notag\\
\lambda_2\leq p_2\notag\\
\lambda_1+\lambda_2\leq p_1 + p_2 -p_1p_2 \label{excapineq}
\end{gather}\vspace{-1cm}

\begin{figure}[h]
    \centering
    \psfrag{a}[][][1]{$\lambda_1$}
    \psfrag{b}[][][1]{$\lambda_2$}
    \psfrag{c}[][][1]{$p_1$}
    \psfrag{d}[][][1]{$p_2$}
    \includegraphics[width=0.35\textwidth]{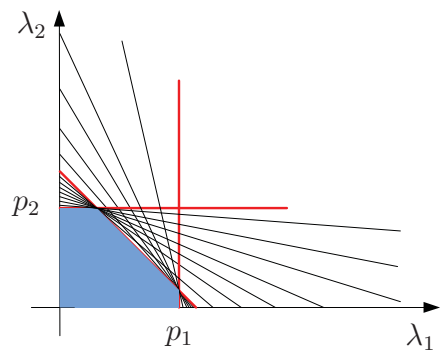}
	\caption{ Stability region for the queueing system of Figure \ref{example}} \vspace{-5mm} 
	\label{capp}
\end{figure}

In Lemmas \ref{l1} to \ref{l4}, we characterized the stability region by an infinite number of its face defining hyperplanes. Figure \ref{capp} depicts the stability region given in (\ref{excapineq}) as well as some of its \textit{face defining hyperplanes} derived according to (\ref{l4_1}) for some sample values of $\alpha=(\alpha_1,\alpha_2)$.  

According to the redundancy theorem in polytopes (refer to Section \ref{polytopesconcepts}), not all the face defining hyperplanes of a polytope are required to describe a polytope. In particular, just the inequalities corresponding to the facets of a polytope are sufficient to characterize a polytope. 
In the following, we will characterize the facets of polytope $\cal P$. The analysis follows and constitutes the main contribution of the paper.

Let $V$ denote the set of vectors $\alpha \in \mathbb{R}^N_+$ with the following property, i.e.,
\begin{eqnarray} \label{Vdef}
\lefteqn{V:=\left\{ \alpha \in \mathbb R^N_+ \mid \forall~ ( \mathcal{U} \subset{\cal N}, \mathcal{U}\neq\emptyset, \alpha_{\mathcal{U}}\neq \underline{0}_{|\mathcal{U}|}, \alpha_{\mathcal{U}^c}\neq \underline{0}_{|\mathcal{U}^c|}) \right.}\nonumber\\&& 
 \left. ~ \exists ~( i \in  \mathcal{U} , j \in  \mathcal{U}^c ,~ m,n \in \mathcal M ,~ \alpha_{i}, \alpha_{j}, m,n\neq 0 ) :  
 \alpha_{i}m=\alpha_{j}n  \right\}.
\end{eqnarray}
A vector $\alpha$ belongs to set $V$ if for any partitioning of elements of vector $\alpha$ into two non-empty disjoint (sub)vectors in which all the elements of each vector are not identically equal to zero, there exists at least one non-zero element in each (sub)vector, the ratio of which is equal to the ratio of two non-zero elements of set $\cal M$. To clarify the definition of set $V$ consider the following example.

\el
 Let ${\cal M}=\{0,1,2\}$ and $N=4$. Now, consider vector $\alpha_1=(1,2, 5 ,10)$. We can partition $\alpha_1$ to $\alpha_{1_{\{1,2\}}}=(1,2)$ and $\alpha_{1_{\{3,4\}}}=(5,10)$. Note that there exist no two elements one in $\alpha_{1_{\{1,2\}}}$ and the other in $\alpha_{1_{\{3,4\}}}$ whose ratio is 1 or 2 and therefore, $\alpha_1 \notin V$. Another example is vector $\alpha_2=(0,2.5,0,0)$ in which for any partitioning of the vector into two non-empty disjoint (sub)vectors one will always have all the elements identically equal to zero. Therefore, $\alpha_2 \in V$.

In the following, we will show that any vector in ${\mathbb R^N_+}-V$ cannot form a facet defining hyperplane for polytope $\cal P$ and therefore the inequalities derived by the vectors of set ${\mathbb R^N_+}-V$ are redundant. This is shown in the following lemma.

\begin{lem} \label{l6}
If $\alpha \in {\mathbb R^N_+}-V$, the hyperplane associated to the valid inequality of (\ref{l4_1}) is not a facet defining hyperplane of $\cal P$ .
\end{lem}
The proof is brought in Appendix \ref{a1}.

Note that $V$ is an infinite set. However, we can show that set $V$ is finite up to multiplication of vectors by positive scalars. To show this, we define set ${\cal W}=\left\{z \in \mathbb{Z_+} \mid  z=\prod_{j=1}^{N-1}m_j,~m_j \in \cal M \right\}$ and $\mathcal{W}^N=\{(\alpha_1,\alpha_2,...,\alpha_N)  \mid \alpha_n \in \mathcal{W}\} $. Then, we can prove the following Lemma.
\begin{lem}
\label{l5}
There exists $\widehat{V} \subseteq {\cal W}^N$ such that any vector $\alpha \in V$ can be written as $\alpha=q \beta$ for some vector $\beta \in \widehat{V}$ and scalar $q>0$.
\end{lem}
The proof is given in Appendix \ref{a2}. 

Recalling the two-queue single-server example of Figure \ref{example}, using Lemmas \ref{l6} and \ref{l5} we can characterize a finite subset of $\mathbb{R}_+^2$ that can produce all the normal vectors of the facet defining hyperplanes of the stability region polytope. For the aforementioned example since $M=1$ we have $\mathcal{W}=\{0,1\}$ and therefore $\mathcal{W}^2=\{(0,0),(0,1),(1,0),(1,1)\}$ and therefore $\widehat{V}\subseteq \{(0,1),(1,0),(1,1)\}$. For this simple example since $M=1$ (as we will see in Corollary \ref{onoff} below), we have $\widehat{V}=\mathcal{W}^N-\{\underline{0}_N\}$. The facet defining hyperplanes of the stability region for this example as described by the vectors in the set $\widehat{V}$ are shown with bold red lines in Figure \ref{capp}.  

\textbf{Fact}: The size of $\cal W$ is $\left|{\cal W}\right|= {M+N-2 \choose{N-1}}+1$. 

In order to prove this fact, note that $\cal M$ has $M$ non-zero elements. Each non-zero element of $\cal W$ comes from the multiplication of $N-1$ non-zero elements of $\cal M$. This is equivalent to the counting problem of choosing $N-1$ balls from $M$ distinctly marked balls without ordering and with replacement which is equal to $\left|{\cal W}\right|= {M+N-2 \choose{N-1}}$ \cite{garcia}. Since $0 \in \cal W$ we must increase this number by 1. 
Hence, $|\mathcal{W}^N|= \left( {M+N-2 \choose{N-1}}+1 \right)^N $. Since $\widehat{V} \subseteq {\cal W}^N -\{\underline{0}_N\} $, therefore 
$|\widehat{V}| \leq |\mathcal{W}^N|-1$ (excluding $\underline{0}_N$ as it results into the obvious equality $0=0$).

In the theorem that follows, we introduce the necessary conditions for the stability of MQMS system with stationary channel distribution.
\begin{theor} 
\label{t1}
 If there exists a server allocation policy under which the system is stable, then 
\begin{eqnarray}
\label {t1_1} \alpha \overline{E[A(t)]^{\mathrm{T}}} \leq \displaystyle\sum_{s\in \mathcal{S}}\pi_s 
  \max_{I \in \mathcal I}  \left(\alpha (C_s \circledast I )\underline{1}_K^{\mathrm{T}} \right),~~ \alpha \in \widehat{V}.
\end{eqnarray}
\end{theor}
\begin{IEEEproof}
The proof follows directly from Lemmas \ref{l4} to \ref{l5}. In fact, from Lemma \ref{l4} we can characterize polytope $\cal P$ by an infinite number of inequalities ($\alpha \in \mathbb{R}_+^N$). In Lemmas \ref{l6} and \ref{l5} we showed that not all of such $\alpha$ vectors are necessary to describe the facets of polytope $\cal P$. More specifically, we proved that just $\alpha \in \widehat{V}$ ($|\widehat{V}| < \infty$) are sufficient to describe the facets of polytope $\cal P$ and other $\alpha$ vectors make redundant inequalities and therefore the theorem follows.
\end{IEEEproof}

According to Theorem \ref{t1}, in order to characterize polytope $\cal P$, we have to specify all the elements of set $\widehat{V}$ which can be obtained after considering all the possible vectors $\alpha \in  \mathcal{W}^N- \{\underline{0}_N\}$ and removing redundancies following (\ref{Vdef}). This can be achieved numerically. 
We may also avoid the numerical computation to derive set $\widehat{V}$ and instead use set $\mathcal{W}^N-\{\underline{0}_N\}$ which is a finite superset of $\widehat{V}$ that contains redundancies (produces redundant inequalities).
%
Although by using $ \mathcal{W}^N- \{\underline{0}_N\} $ we may obtain some redundant inequalities, since $|\mathcal{W}^N- \{ \underline{0}_N\}|< \infty$ we still have a finite number of inequalities to describe polytope $\mathcal{P}$. 
Table \ref{table} depicts $\widehat{V}$ and $|\mathcal{W}^N|-1$ for some sample cases $N=2,3$ and $M=1,2,3,4$. \vspace{-2mm}
\begin{table}[h]
\renewcommand{\arraystretch}{1.3}
\caption{$|\widehat{V}|$ and $|\mathcal{W}^N|-1$ for $N=2,3$ and $M=1,2,3,4$}
\label{table}
\centering 
\begin{tabular}{cc|c|c|c|c|l}
\cline{3-6} 
& & $M=1$ & $M=2$ & $M=3$ & $M=4$ \\ \hline
\cline{1-6} 
\multicolumn{1}{|c||}{\multirow{2}{*}{\hspace{-1mm}$|\widehat{V}|,|\mathcal{W}^N|-1 $\hspace{-1mm}}} &
\multicolumn{1}{|c|} {$N=2$} & $3,3$ & $5,8$ & $9,15$ &$12,24$     \\ \cline{2-6}
\multicolumn{1}{|c||}{}                        &
\multicolumn{1}{|c|}{$N=3$} & \hspace{-.5mm}$7,7$\hspace{-.5mm} & $25,63$ & $109,342$ & \hspace{-.5mm}$253,1330$\hspace{-.5mm}  \\
%
\hline
\end{tabular}
\end{table}

\begin{cor} \label{nec-stat}
 For an MQMS system with stationary arrival processes, we have $E[A(t)]=\lambda=(\lambda_1,\lambda_2,...,\lambda_N)$ and the necessary conditions for the stability of the system would be
\begin{eqnarray}
\label {t1_2}   \alpha {\lambda}^{\mathrm{T}} \leq \displaystyle\sum_{s\in \mathcal{S}}\pi_s 
  \max_{I \in \mathcal I}  \left(\alpha (C_s \circledast I )\underline{1}_K^{\mathrm{T}} \right),~~ \alpha \in \widehat{V}.
\end{eqnarray}
\end{cor}


\begin{cor} \label{onoff}
 For an MQMS system with ON-OFF channels, we have ${\cal W}=\{0,1\}$ and therefore the necessary conditions for the stability of the system would be
\begin{eqnarray}
\label {t1_3} \alpha \overline{E[A(t)]^{\mathrm{T}}} \leq \displaystyle\sum_{s\in \mathcal{S}}\pi_s 
  \max_{I \in \mathcal I}  \left(\alpha (C_s \circledast I )\underline{1}_K^{\mathrm{T}} \right),~~ \alpha \in \{0,1\}^N-\{\underline{0}_N\}.
\end{eqnarray}
In a special case considered in \cite{khodam} where the channels are modeled by independent Bernoulli random variables with $ E[C_{n,k}(t)]=p_{n,k}$, the necessary conditions for the stability of the system are given by \vspace{-3mm}
 \begin{eqnarray}
\label {t1_4} 
\overline{\displaystyle\sum_{n\in Q}E[A_n(t)]}\leq K-\displaystyle\sum_{k=1}^{K}\displaystyle\prod_{n\in Q}(1-p_{n,k})~~~ \forall Q\subseteq \mathcal{N}.
\end{eqnarray}
 In this case, the total number of inequalities needed to describe polytope $\cal P$ is equal to $2^N-1$.
\end{cor}

\subsection{Sufficient Conditions for the Stability of MQMS System}
 
Consider a server allocation policy which determines the allocation matrix at each time slot $t$ by solving the following maximization problem.
\begin{eqnarray}
\label{mw}
I(t) = \arg\max_{I \in \cal I} X(t-1) (C(t) \circledast I )\underline{1}_K^{\mathrm{T}}
\end{eqnarray}
This policy is called Maximum Weight (MW) and was introduced in \cite{javidi2,Tassiulas92}. According to the constraints on allocation matrix $I$, we can easily conclude that MW policy allocates the servers by the following rule:
At each time slot $t$, each server $k$ is allocated to the queue $n$ that achieves the maximum $X_n(t-1)C_{n,k}(t)$. 

In the special case where the channels are ON-OFF, the policy is the same as AS/LCQ (Any Server/Longest Connected Queue) introduced in \cite{khodam}. In the following, we derive the sufficient conditions for the stability of our model and prove that MW stabilizes the system as long as condition (\ref{sc1}) below is satisfied. We also derive an upper bound for the time averaged expected number of packets in the system.
\begin{theor}\label{sc} The MQMS system is stable under MW policy if \emph{for all} \(t\)
\begin{equation}
\label{sc1}\alpha E[A(t)]^{\mathrm{T}} < \displaystyle\sum_{s\in \mathcal{S}}\pi_s 
  \max_{I \in \mathcal I}  \left(\alpha (C_s \circledast I )\underline{1}_K^{\mathrm{T}} \right),~~ \alpha \in \widehat{V}.
\end{equation}
Furthermore, the following bound for the average expected ``aggregate" occupancy holds.
\begin{eqnarray}
\label{sc2} \limsup_{t \to \infty} {1\over{t}} \displaystyle\sum_ {\tau =0}^{t-1} \displaystyle\sum_{n=1}^N E[X_n(\tau)]\leq \frac{N A_{max}^2+(MK)^2}{2\delta} 
\end{eqnarray}
In (\ref{sc2}), $\delta$ is the maximum positive number such that \textit{for all t} we have $E[A(t)]+\delta \underline{1}_N \in  \mathcal{P}$.
\end{theor}
The detailed proof is given in Appendix \ref{a3}.

\cor \label{suf-stat}
 For an MQMS system with stationary arrival processes, we have $E[A(t)]=\lambda$ and the sufficient conditions for the stability of the system under MW policy are $\alpha {\lambda}^{\mathrm{T}} < \sum_{s\in \mathcal{S}}\pi_s  \max_{I \in \mathcal I}  \left(\alpha (C_s \circledast I )\underline{1}_K^{\mathrm{T}} \right),~~ \alpha \in \widehat{V}$.

\begin{figure}[tp]
\centering
\def\tabularxcolumn#1{m{#1}}
\begin{tabular}{cc}
\subfloat[$N=2$, $M=1$]{\includegraphics[width=0.35\textwidth]{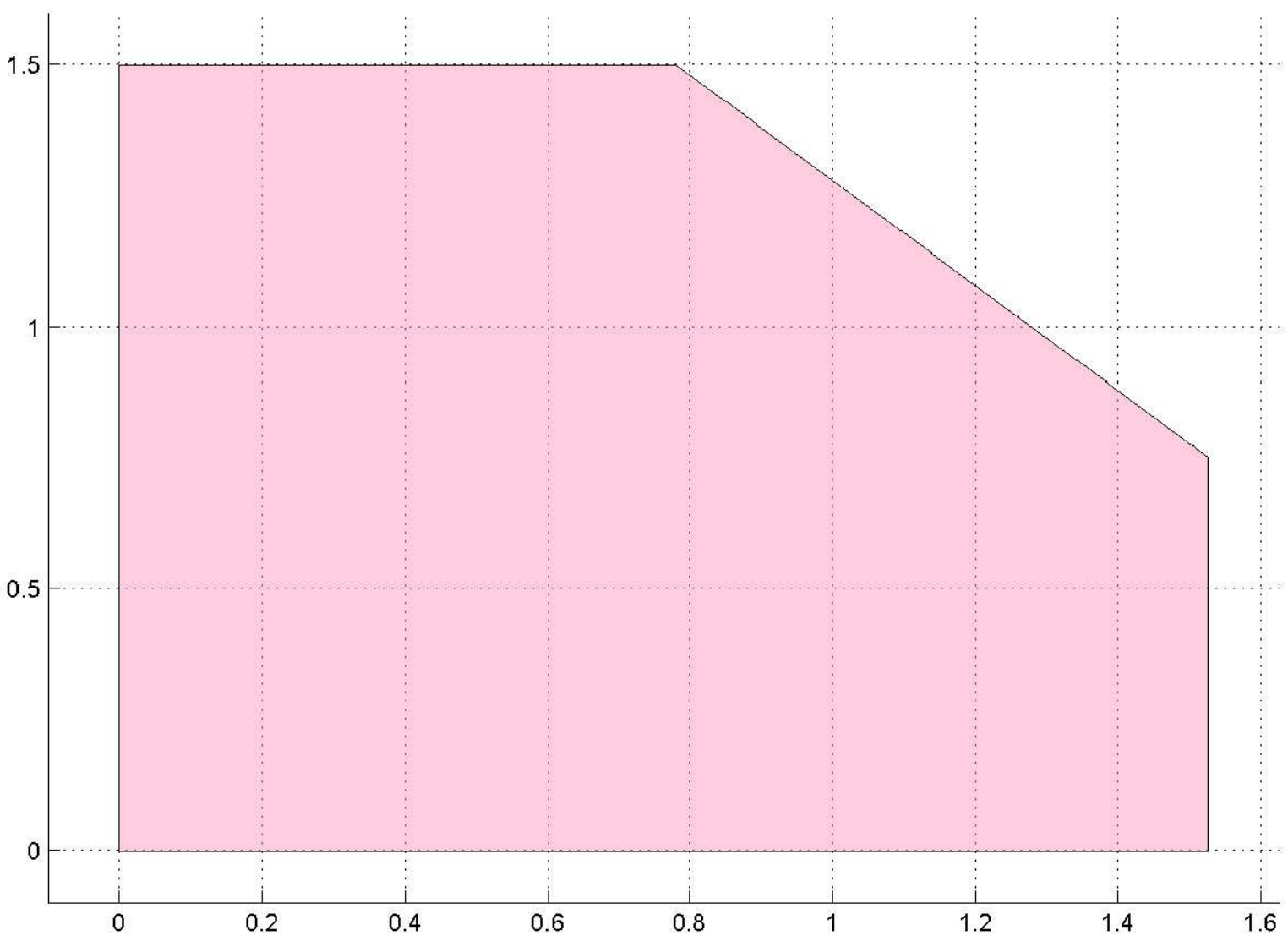}} 
   &\hspace{5mm} \subfloat[$N=3$, $M=1$]{\includegraphics[width=0.35\textwidth]{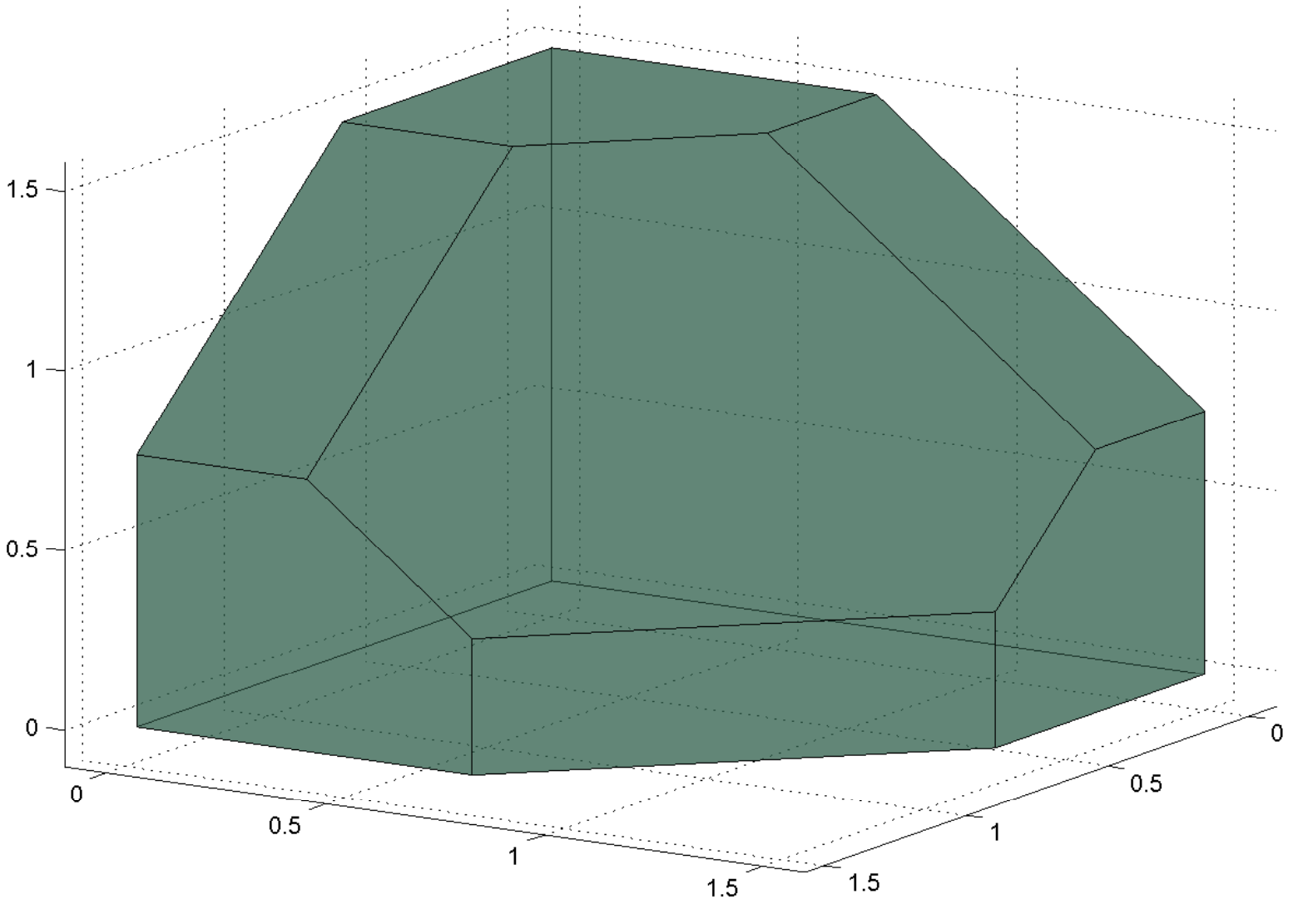}}\\
\subfloat[$N=2$, $M=2$]{\includegraphics[width=0.35\textwidth]{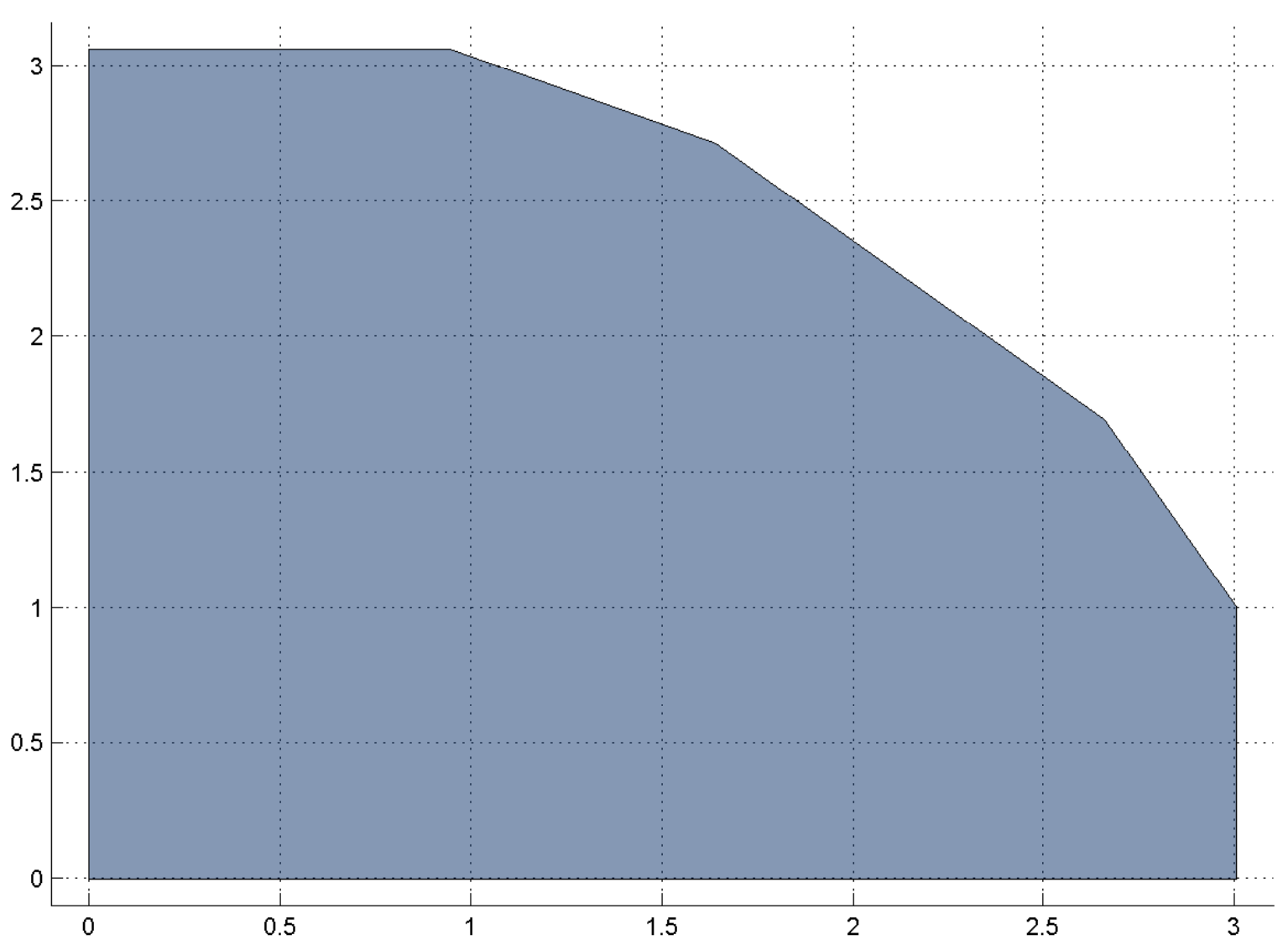}} 
   &\hspace{5mm} \subfloat[$N=3$, $M=2$]{\includegraphics[width=0.35\textwidth]{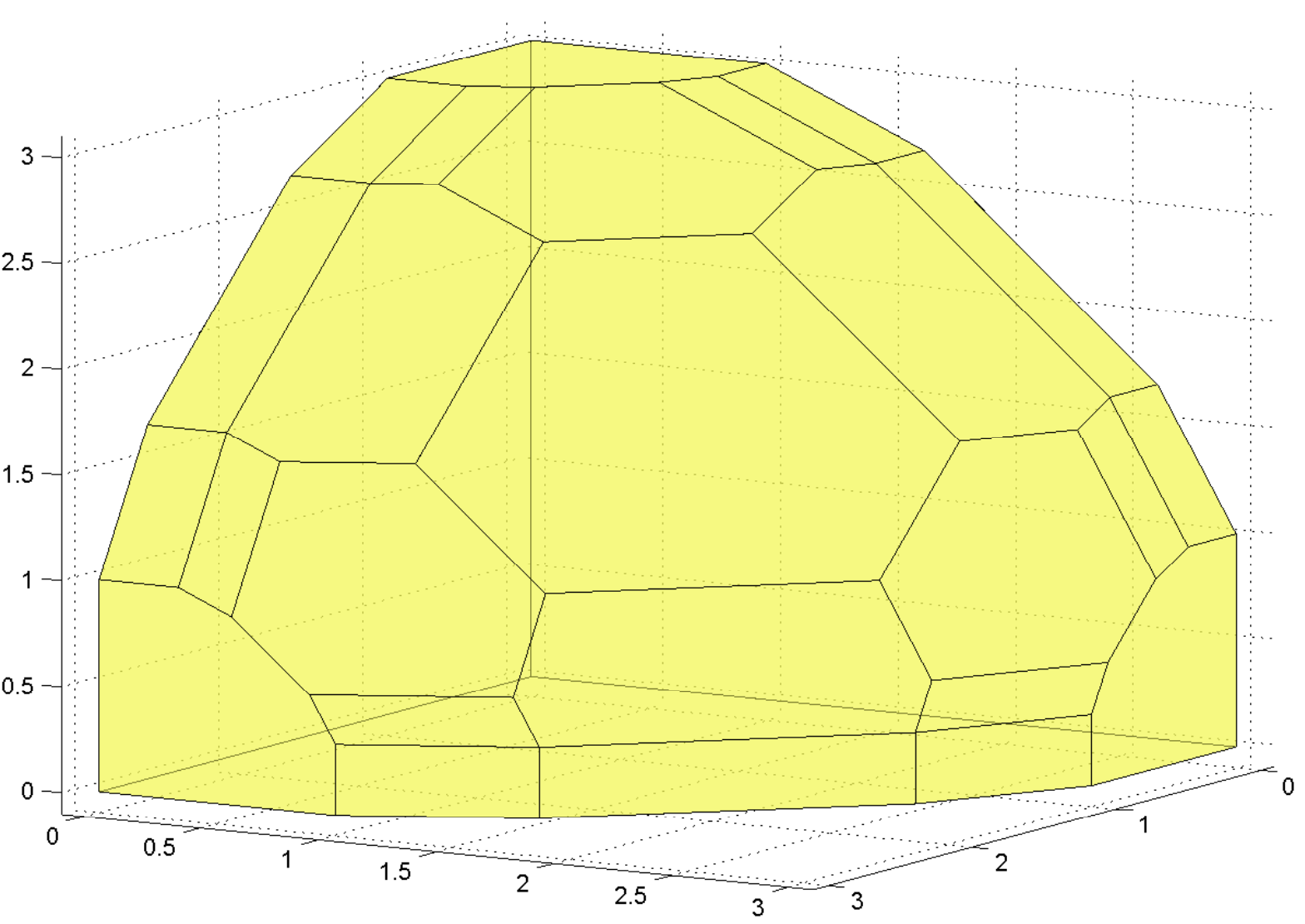}}\\
\subfloat[$N=2$, $M=3$]{\includegraphics[width=0.35\textwidth]{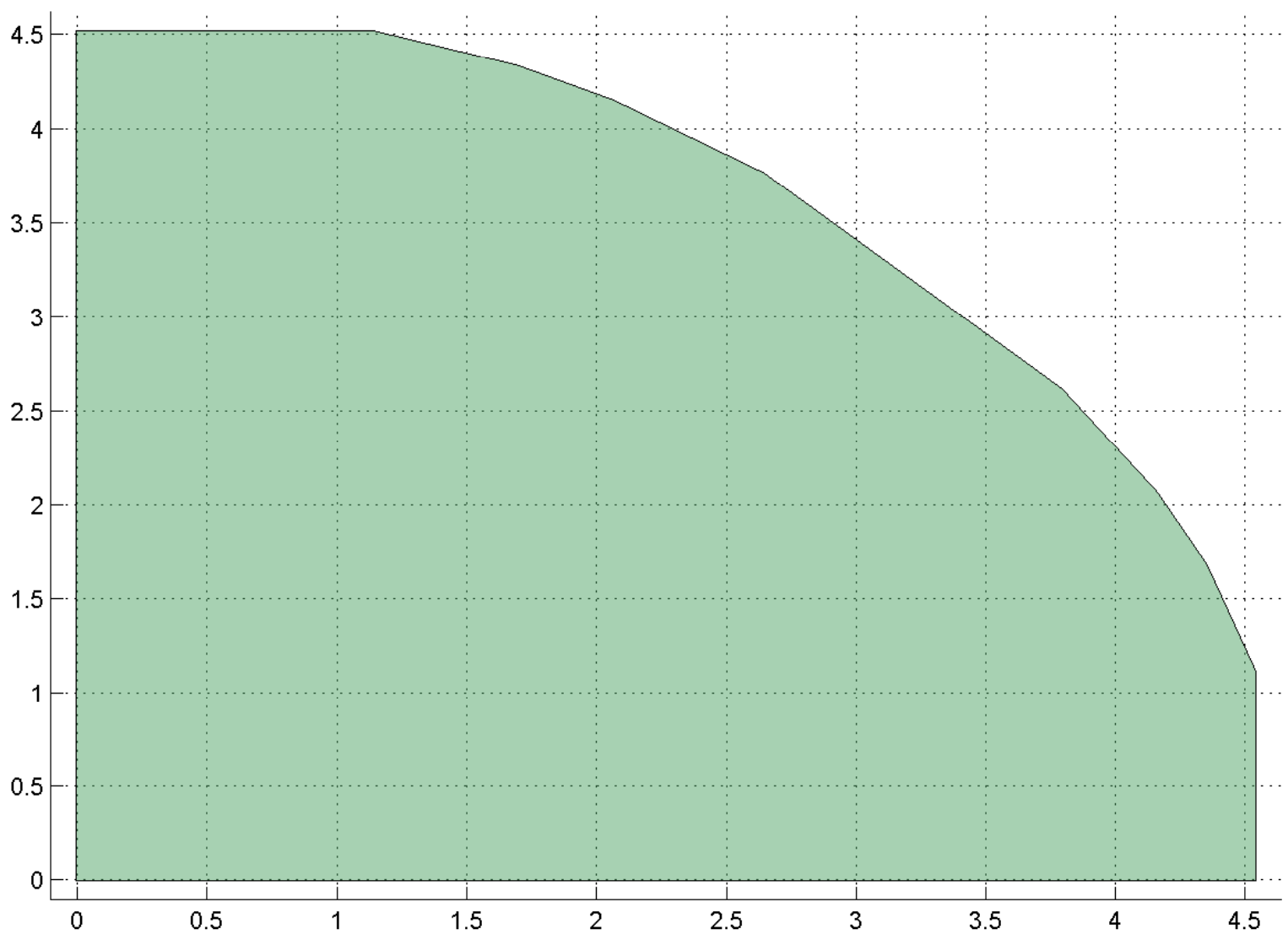}} 
   &\hspace{5mm} \subfloat[$N=3$, $M=3$]{\includegraphics[width=0.35\textwidth]{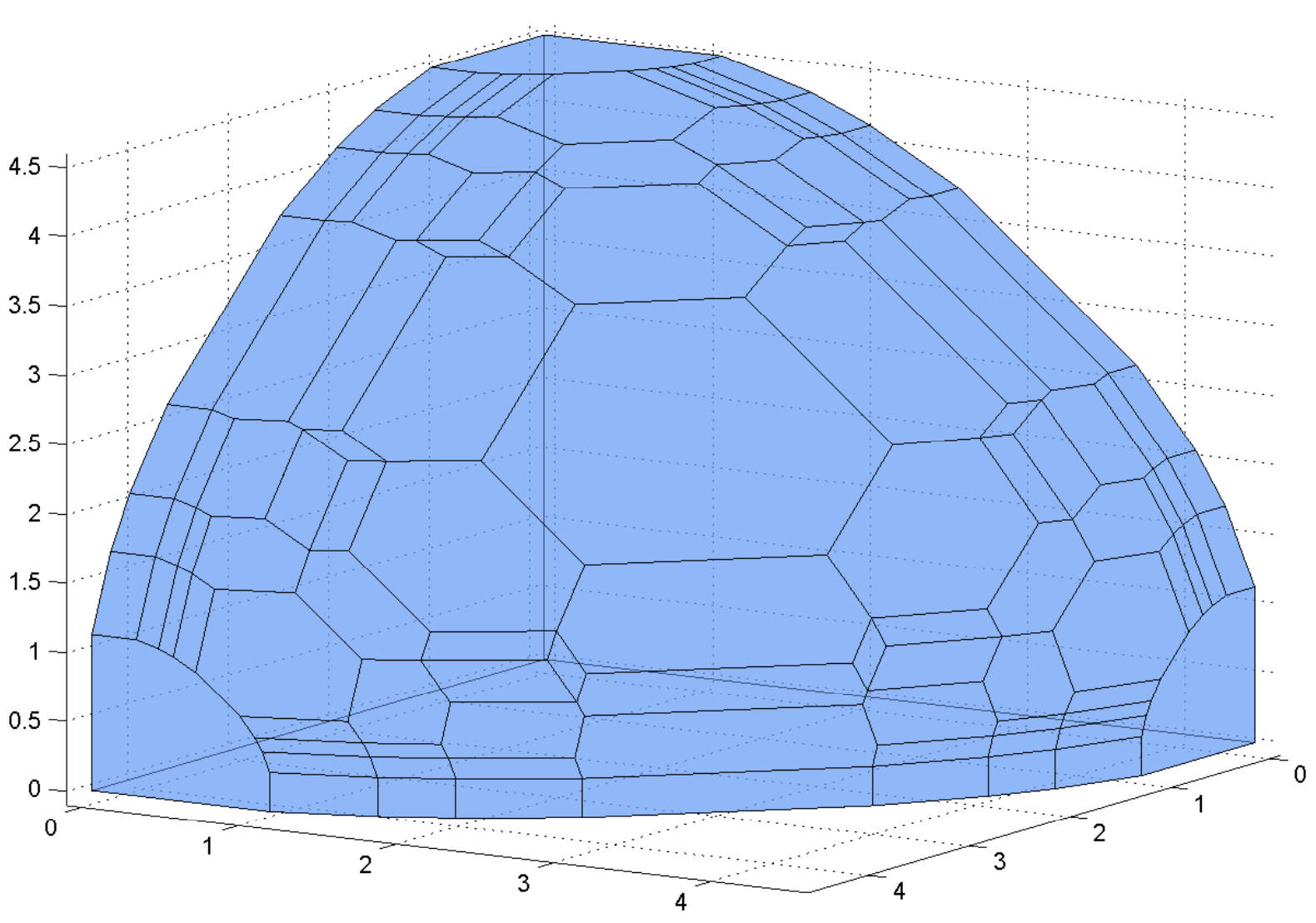}}\\
\subfloat[$N=2$, $M=4$]{\includegraphics[width=0.35\textwidth]{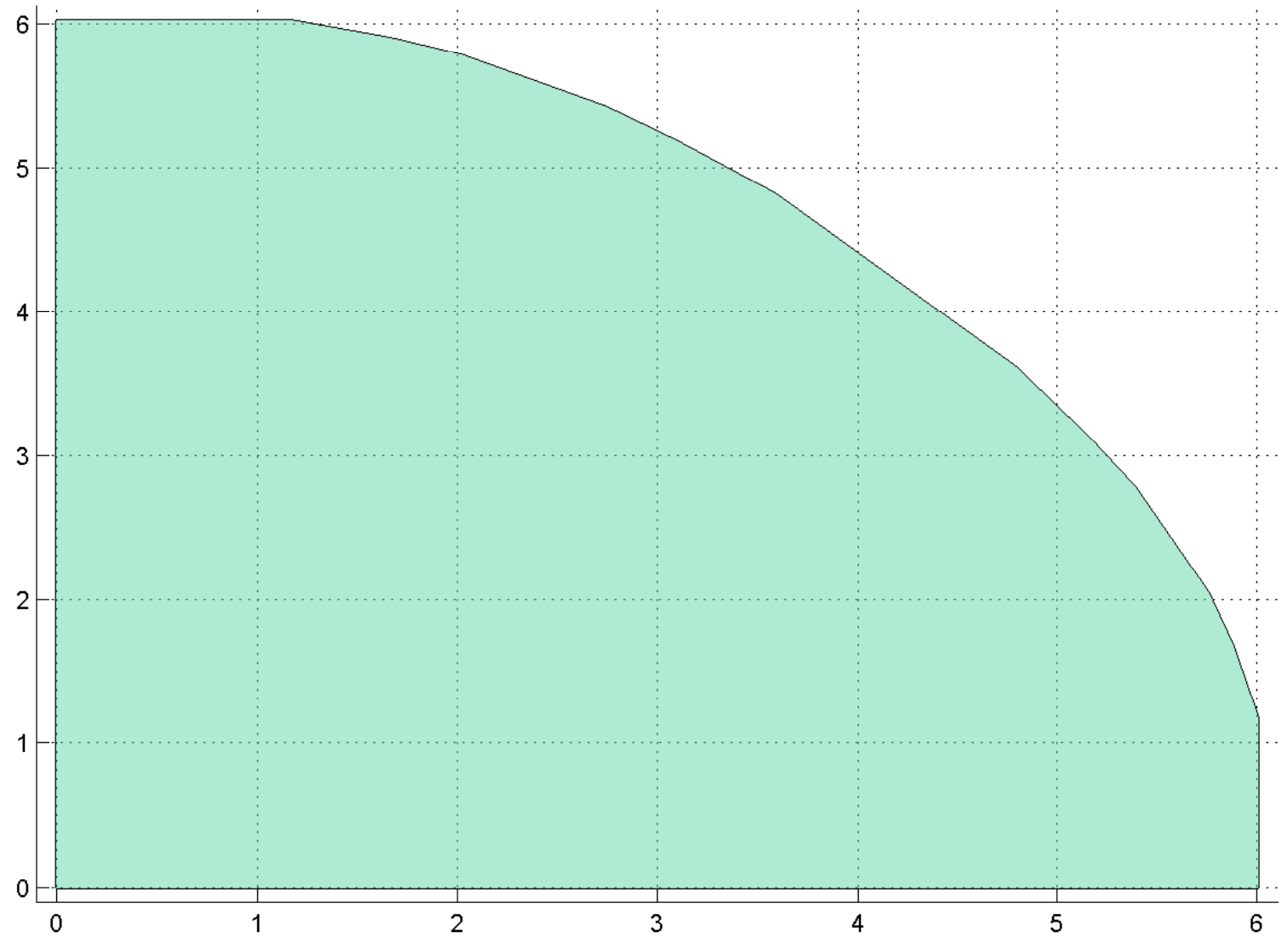}} 
   &\hspace{5mm} \subfloat[$N=3$, $M=4$]{\includegraphics[width=0.35\textwidth]{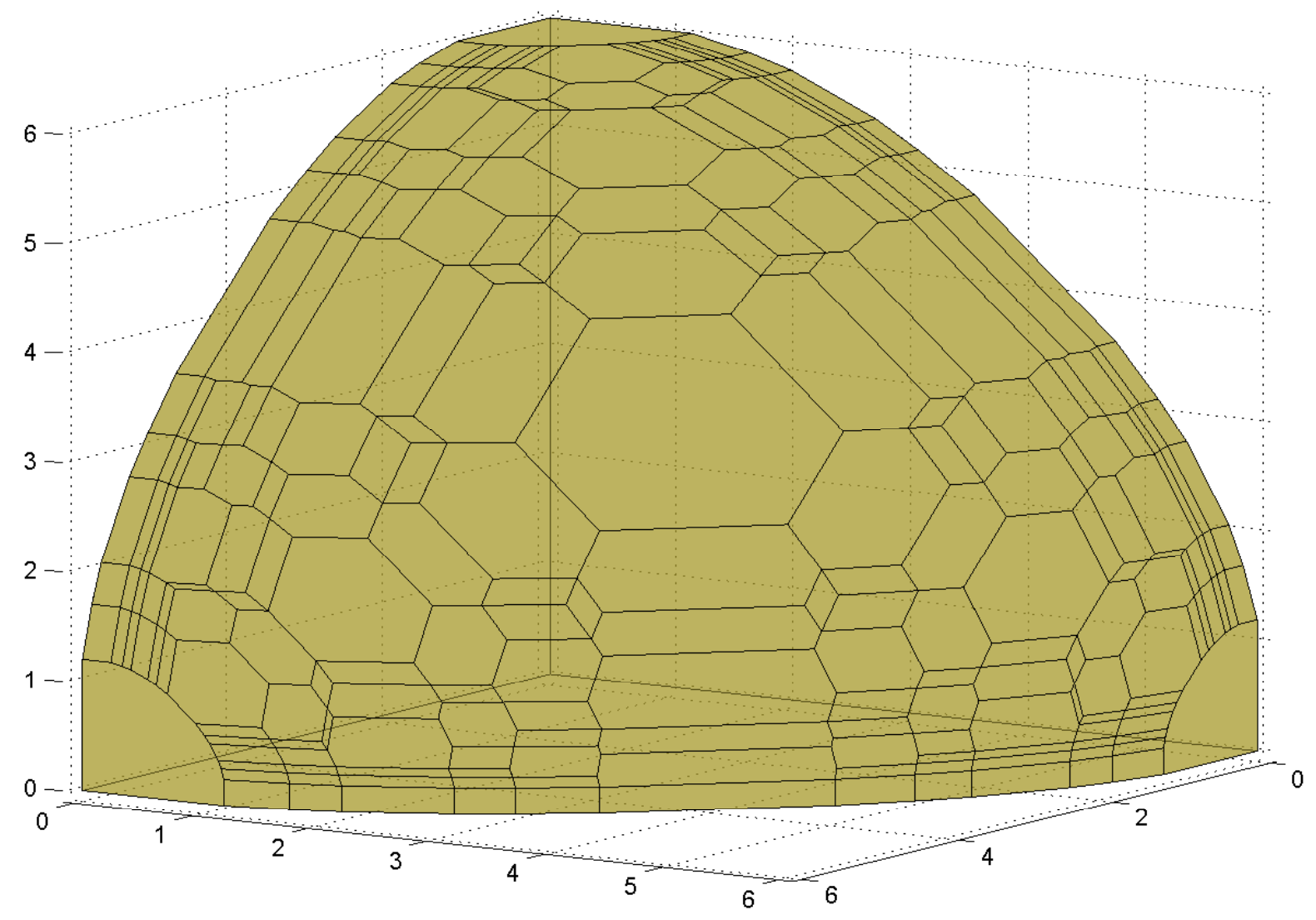}}\\
\end{tabular}
\caption{Stability region for $N=2,3$, $M=1,2,3,4$ and $K=3$ }\label{cap}
\end{figure}

According to Corollaries \ref{nec-stat} and \ref{suf-stat} and the definition of network stability region, we can conclude that for an MQMS system with stationary channel distribution and stationary arrival processes, the stability region is characterized by (\ref{t1_2}). The stability region described in (\ref{t1_2}) is a polytope which can be imagined for $2$ and $3$ dimensional systems, i.e., $ N=2,3 $. Figure \ref{cap} shows the stability region for $ N=2,3$ and $M=1,2,3,4$. In all the cases $K=3$.

Consider an MQMS system with ON-OFF channels. For such a system, MW policy is equivalent to AS/LCQ policy introduced in \cite{khodam}. AS/LCQ policy takes an arbitrary ordering of servers (to be allocated to the queues) and then for each server allocates it to its longest connected queue (LCQ). 
Note that in such a system, ${\cal W}=\{0,1\}$ and therefore we have the following corollary.

\begin{cor}
The MQMS system with ON-OFF channels is stable under AS/LCQ if \textit{for all} $t$
\begin{eqnarray}
\label {t2_3} 
 \alpha E[A(t)]^{\mathrm{T}} < \displaystyle\sum_{s\in \mathcal{S}}\pi_s 
  \max_{I \in \mathcal I}  \left(\alpha (C_s \circledast I )\underline{1}_K^{\mathrm{T}} \right),~~ \alpha \in \{0,1\}^N-\{\underline{0}_N\}. \vspace{-1mm}
\end{eqnarray}
In a special case where the channels are modeled by independent Bernoulli random variables with $ E[C_{n,k}(t)]=p_{n,k}$, AS/LCQ stabilizes the system if \textit{for all} $t$\vspace{-2mm}
\begin{eqnarray}
\label {t2_4} 
\displaystyle\sum_{n\in Q}E[A_n(t)] < K-\displaystyle\sum_{k=1}^{K}\displaystyle\prod_{n\in Q}(1-p_{n,k})~~~ \forall Q\subseteq \mathcal{N}.
\end{eqnarray}
\end{cor}

\subsection{Numerical Example}
In this section, we will consider a simple numerical example of MQMS system with $N=2$, $K=3$ and $M=3$.  For such as system the size of channel state space is $4^6=4096$. We have chosen the channel distribution randomly. We also assume that the utility functions associated to queues $1$ and $2$ are $\log(1+10r_1)$ and $10r_2$, respectively\footnote{Linear and logarithmic utility functions were studied in \cite{Now,phd}.}. Therefore the total utility function is $f(r_1,r_2)=\log(1+10r_1)+10r_2$. We also assume that the packet arrival processes to queues $1$ and $2$ are Poisson distributed with rate $5$ packets/time slot. First, we are going to characterize the stability region for this example using the theorems given in the previous subsections. Then, we will correlate our findings with the results from \cite{phd,Now,neely-infocom05}.

\begin{figure}[h]
    \centering
    \includegraphics[width=0.65\textwidth]{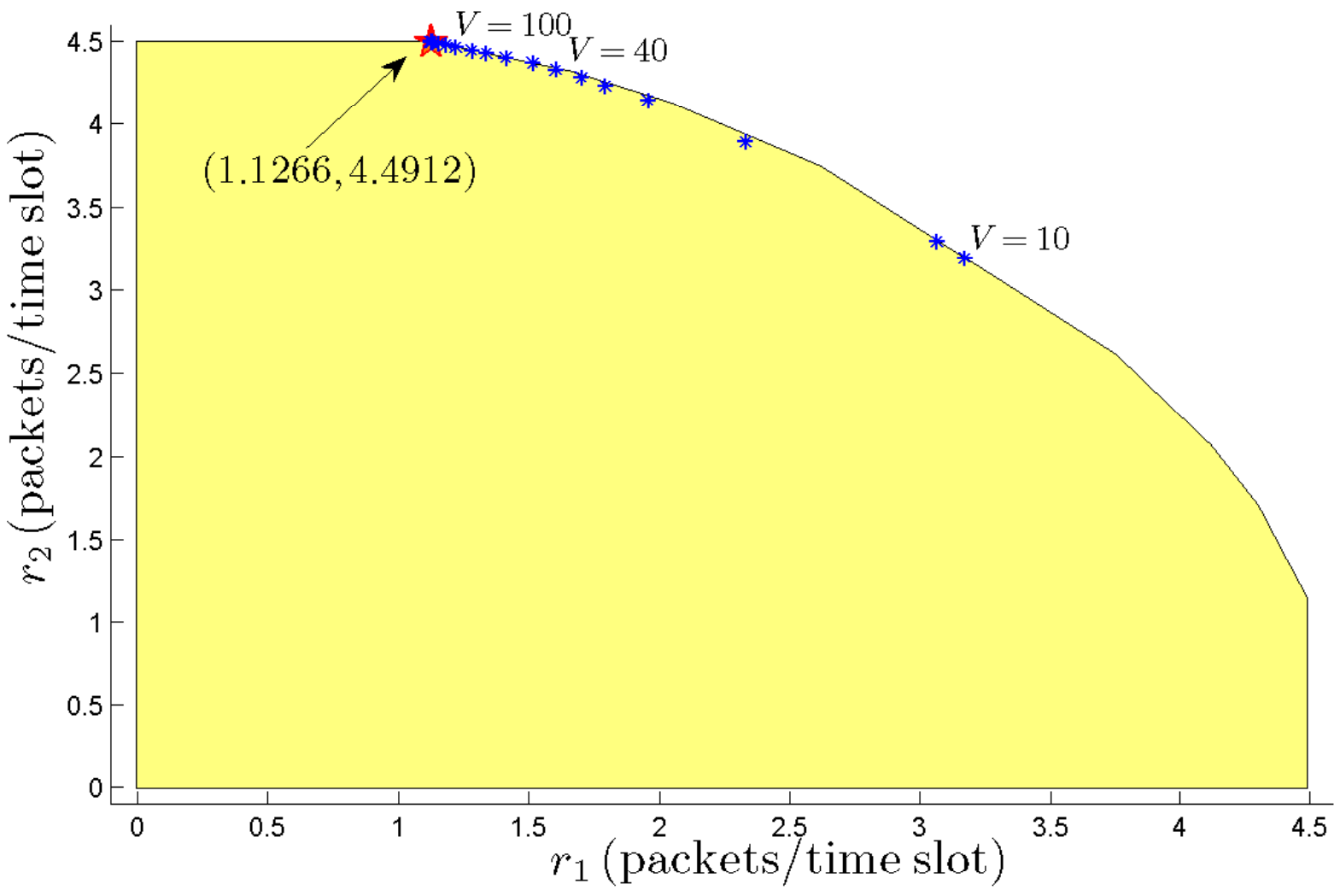}
	\caption{ Stability region for $N=2$, $K=3$, $M=3$ with a random channel distribution (the optimal solution and the solution of CLC2b algorithm are also depicted in this figure)}
	\label{ex-cap} \vspace{-7mm}
\end{figure}

For the described system we can easily check that $\mathcal{W}=\{0,1,2,3\}$. Thus, $\mathcal{W}^2-\{(0,0)\}=\{(0,1),(0,2),(0,3),(1,0),(1,1),(1,2),(1,3),(2,0),(2,1),(2,2),(2,3),(3,0),(3,1),(3,2),(3,3)\}$ whose size is $15$ (as mentioned in Table \ref{table}). Note that some of these vectors make redundant inequalities. For example using $(1,1)$, $(2,2)$ and $(3,3)$ as the vector $\alpha$ in (\ref{t1_2}) will result in the same inequalities. By removing the redundant vectors we obtain set $\widehat{V}=\{(0,1),(1,0),(1,1),(1,2),$ $(1,3),(2,1),(2,3),(3,1),(3,2)\}.$ 
Using (\ref{t1_2}), the stability region is characterized by the following set of inequalities which is also shown in Figure \ref{ex-cap}.
\begin{align} \label{constr}
r_1 \leq 4.4792 &&  r_1+2r_2 \leq 10.2893 && 2r_1+3r_2 \leq 16.4564 \notag\\
r_2 \leq 4.4912 &&  r_1+3r_2 \leq 14.6002 && 3r_1 + r_2 \leq 14.5803 \notag \\
r_1+r_2 \leq 6.3577 && 2r_1+r_2 \leq 10.2874 && 3r_1+ 2 r_2 \leq 16.4639
\end{align}
Therefore, the utility optimization problem is specified as 
\begin{eqnarray}\label{opt-example}
\label{util}~~ \texttt{Maximize:}~~~~~~~~~~ \log(1+10r_1)+10r_2~~~~~~~~~~~~~~~~~~~~~~~~~~~~~~~~~~~~~~~~~~~~~~~~~~~~~~  \\ 
\texttt{Subject to:}~~~~~~~~\text{All the constraints in (\ref{constr}})~~~~~~~~~~~~~~~~~~~~~~~~~~~~~~~~~~~~~~~~~~~~~~~~ \notag \\
0 \leq r_1 \leq 5~~~,~~~ 0 \leq r_2 \leq 5.~~~~~~~~~~~~~~~~~~~~~~~~~~~~~~~~~~~~~~~~~~~~~  \nonumber  
\end{eqnarray}
The solution to the above problem is $r^\star=(r_1^\star,r_2^\star)=(1.1266,4.4912)$ and is indicated by the star in bold in Figure \ref{ex-cap}. As we expect the optimal point is located on the boundary of the stability region. For each queue, we can use a leaky bucket and adjust the admitted rate to each queue to the optimal value we obtained from the optimization problem.

We have also simulated the flow control strategy CLC2b (Cross-Layer Control 2-b) which was introduced in \cite{phd,Now,neely-infocom05} (Due to lack of space the reader is suggested to refer to the literature regarding this algorithm). CLC2b was proposed as a flow control strategy for a general queueing system for which the \textit{stability region is unknown}. The algorithm assumes that $A_n(t)\leq R_n^{max}$ for all $t$.
Assume that the backlog in the transport layer of queue $n$ at the beginning of time slot $t$ is denoted by $L_n(t)$.
CLC2b algorithm is a dynamic flow control algorithm in the transport layer in which at each time slot it determines how many packets to admit to queue $n$ (in the network layer) from the transport layer (i.e., from $L_n(t)$ packets in the transport layer waiting for admission). Assume that this number is denoted by $r_n(t)$. Therefore, the goal of CLC2b algorithm is to determine $r_n(t)$ such that $f(\bar{r})$ is maximized while the queues are kept stable. 
CLC2b algorithm incorporates two novel notions, namely \textit{auxiliary variables} and \textit{virtual cost queues}. In the following, we will briefly review them.

\begin{figure}[tp]
    \centering
    \psfrag{a}[][][1]{$L_n(t)$}
    \psfrag{b}[][][1]{$Y_n(t)$}
    \psfrag{c}[][][1]{$X_n(t)$}
    \psfrag{d}[][][1]{$A_n(t)$}
    \psfrag{e}[][][1]{$r_n(t)$}
    \psfrag{f}[][][1]{Transport Layer}
    \psfrag{g}[][][1]{Network Layer}   
    \psfrag{j}[][][1]{$\gamma_n(t)$}
    \includegraphics[width=0.6\textwidth]{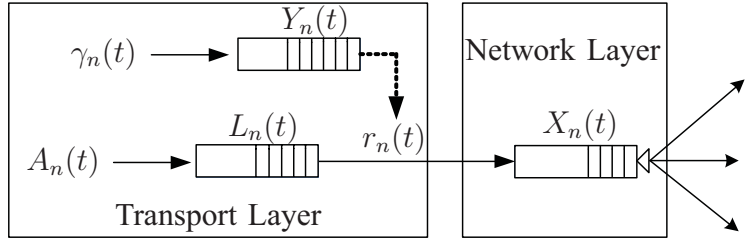}
	\caption{CLC2b algorithm framework}
	\label{flowcontrol} \vspace{-5mm}
\end{figure}

The utility optimization in (\ref{fairness}) can be transformed into the following stochastic optimization problem by introducing the auxiliary variables $\gamma_n, n\in \cal N$ \cite{phd,Now,neely-infocom05}.
\begin{eqnarray}
\label{fairness2}~~ \texttt{Maximize:}_{\substack{\\  \\ \hspace{-2.2cm} \gamma , r}}~~~~~~~~~~\sum_{n=1}^N f_n(\gamma_n)~~~~~~~~~~~~~~~~~~~~  \\ 
\texttt{Subject to:}~~~~~~~~~~r_n\geq \gamma_n~~~~~~~~~~~~~~~~~~~~~~~ \notag\\
r=(r_1,r_2,...,r_N) \in \Lambda ~~~~~\nonumber \\
0 \leq r_n \leq \lambda_n~~\forall n=1,2,...,N. \nonumber\hspace{-.5cm}
\end{eqnarray}
For each queue $n$ a virtual cost queue with virtual arrival process $\gamma_n(t)$ and virtual service process $r_n(t)$ is defined (See Figure \ref{flowcontrol} and refer to \cite{Now}). Note that the service process of the virtual queue $n$ is exactly equal to the admitted packet process to the actual queue $n$. Assuming that the length of this virtual queue is denoted by $Y_n(t)$ at time slot $t$, the following equation reveals the evolution of this process by time.
\begin{eqnarray} \label{yevol}
Y_n(t)=(Y_n(t-1)-r_n(t))^+ +\gamma_n(t)
\end{eqnarray}
It was shown in \cite{Now,phd} that if a control strategy stabilizes both the actual and virtual queues then the resulting averages $\overline{r}_n$ and $\overline{\gamma}_n$ will satisfy the constraints of the optimization problem (\ref{fairness2}).
The flow control strategy of CLC2b is specified as follows.
\begin{itemize}
\item Every timeslot and for each queue $n$ observe $X_n(t-1)$ and $Y_n(t-1)$ and determine 
\begin{eqnarray}
r_n(t)=\left\{\begin{matrix}
\min\{L_n(t)+A_n(t),R_n^{max}\}&&\text{if}~\eta Y_n(t-1)>X_n(t-1)
\\ 
0&&\text{otherwise}
\end{matrix}\right.
\end{eqnarray}
\item Choose $\gamma_n(t)$ (the input process of queue $n$'s virtual queue) as the solution of the following optimization problem and update $Y_n(t)$ following (\ref{yevol}). \vspace{-3mm}
\begin{eqnarray}
\texttt{\begin{small}Maximize:\end{small}}_{\substack{\\ \\ \hspace{-2cm}\gamma}}~~~~~~ Vf_n(\gamma)-\eta Y_n(t-1)\gamma  \\ 
\texttt{\begin{small}Subject to:\end{small}}~~~~~~0 \leq \gamma \leq R_n^{max}.~~~~~~~~\nonumber  
\end{eqnarray}
\item Use the MW policy for resource (server) allocation and update $X_n(t)$.
\end{itemize}
$f_n(\cdot)$ is the utility function associated to queue $n$ and $\eta$ is a positive weight that satisfies $0<\eta \leq 1$. It determines the relative weight of the virtual queue in stabilizing the system. $V$ is a control parameter that affects the proximity of $\overline{r}_n$ (the solution of CLC2b algorithm) to the optimal point $r^\star$. 
It was shown that CLC2b can stabilize both the actual and virtual queues while guaranteeing a lower bound for the achieved utility function (given by (\ref{lbound})) and also an upper bound for the system backlog in the network layer (given by (\ref{ubound})).
%
\begin{eqnarray}\label{lbound}
\liminf_{t \to \infty} f(\overline{r})\geq f(r^\star)-\frac{D}{V}
\end{eqnarray}
\begin{eqnarray} \label{ubound}
\limsup_{t \to \infty} \frac{1}{t}\sum_{\tau=0}^t\sum_{n=1}^{N}X_n(\tau)\leq \frac{D+VG}{\mu}
\end{eqnarray}
Where $D$, $G$ and $\mu$ are positive constants and are dependent to the statistical properties of the system. We can observe that as the control parameter $V$ increases the difference between $f(\overline{r})$ and the optimal utility $f(r^\star)$ decreases as $\frac{1}{V}$. However, the closer $f(\overline{r})$ and $f(r^\star)$ are the larger the backlog is in the system!

\begin{figure}[h]
  \centering
  \subfloat[Difference of the average queue lengths for the optimal solution and the CLC2b algorithm]{\label{queues}\includegraphics[width=0.48\textwidth]{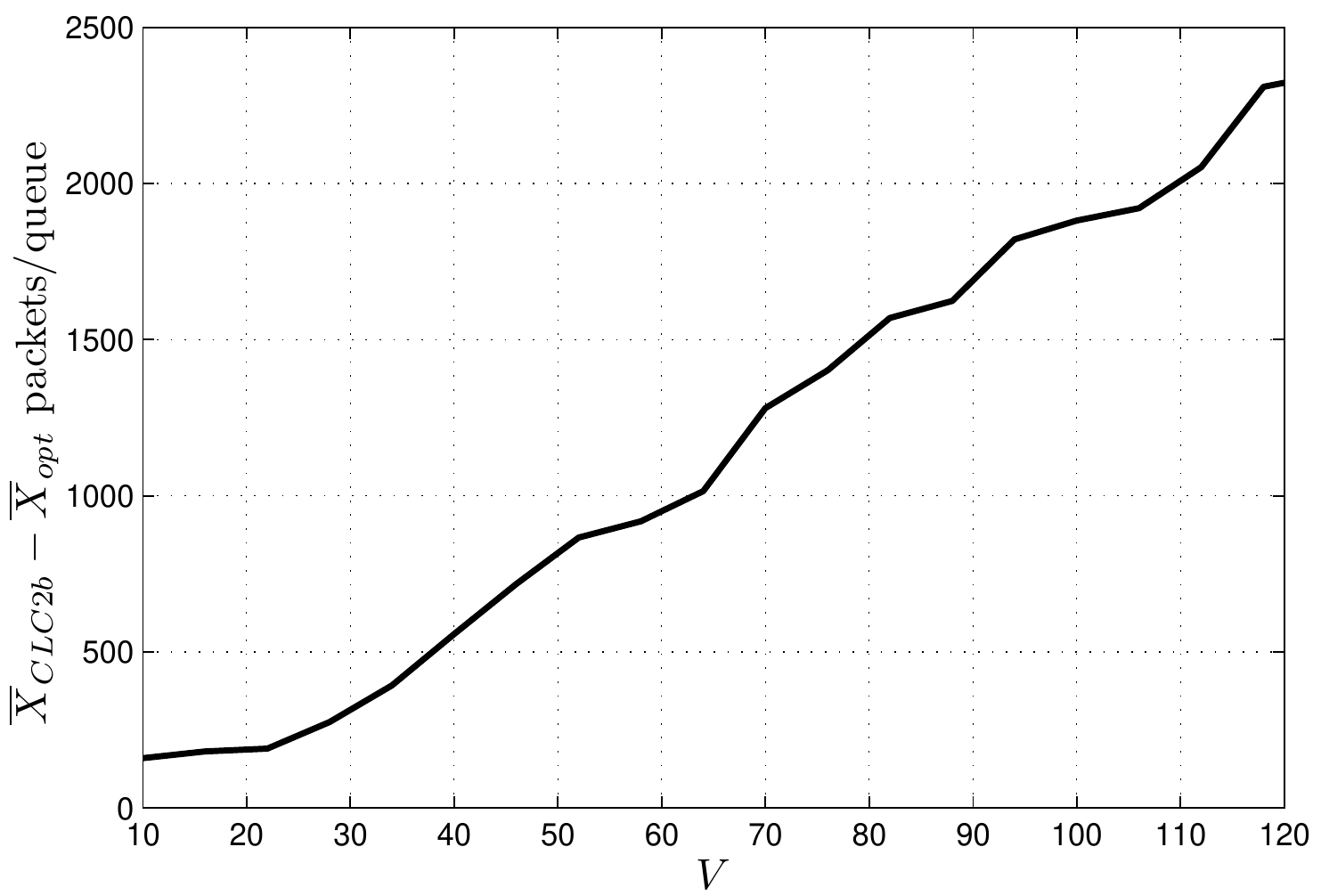}}   \hspace*{5mm}             
  \subfloat[Percentage of the normalized difference of the utility function for $\overline{r}$ and $r^\star$]{\label{util}\includegraphics[width=0.48\textwidth]{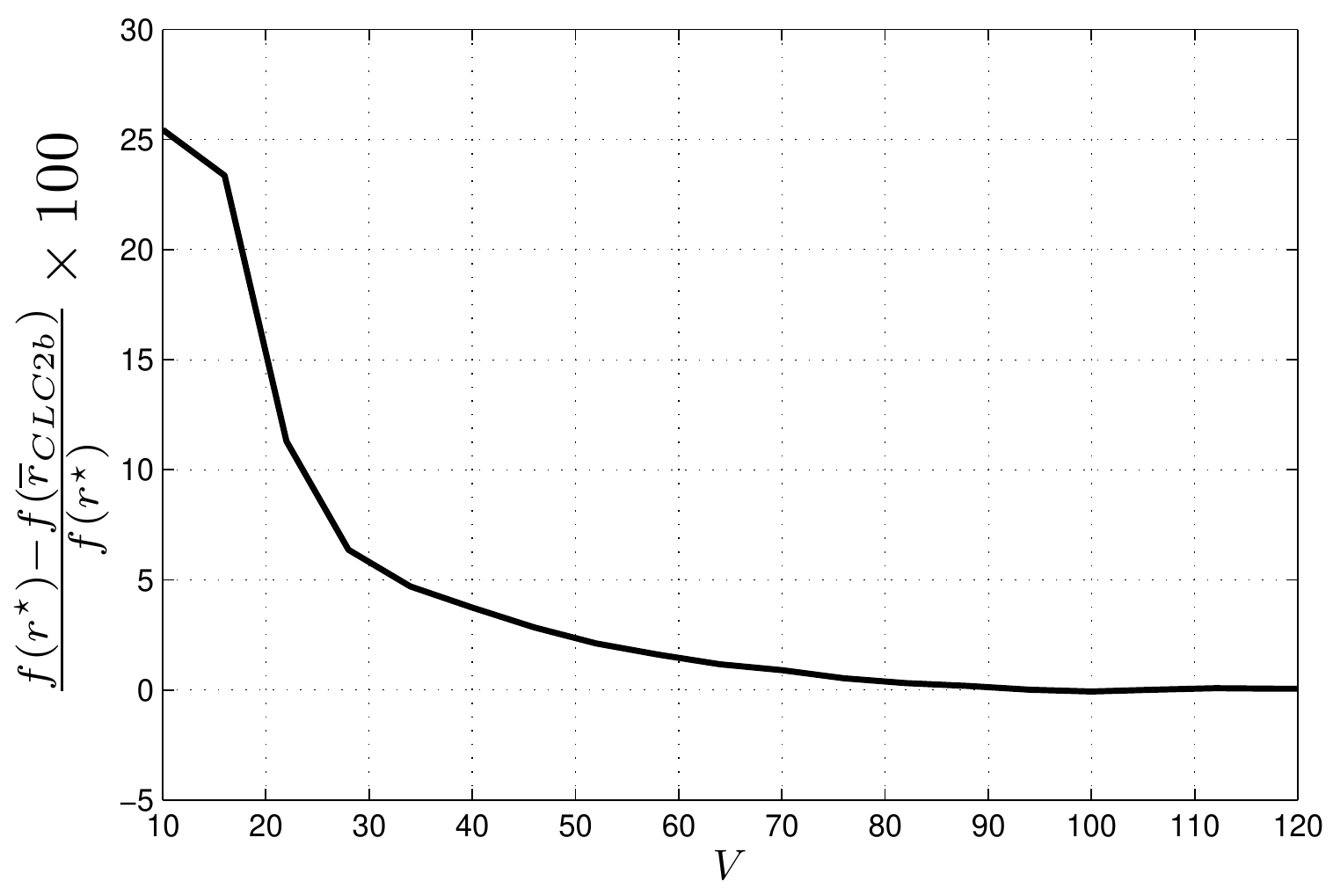}}
  \caption{Comparison of the performances of CLC2b algorithm with the optimal solution} \vspace{-1cm}
\end{figure}

We have simulated the CLC2b algorithm with $\eta=1$. Furthermore, we performed a simulation where the admitted rates to the queues are controlled by leaky buckets outputting the optimal rates $r^\star=(1.1266,4.4912)$ as determined by our analysis earlier.
We compared the two systems by computing the following measures versus $V$.
\begin{itemize}
\item Difference of the average queue occupancy of the systems, i.e., $\overline{X_{CLC2b}}-\overline{X_{opt}}$ shown in Figure \ref{queues}.
\item Percentage of the normalized difference of the utility function for $\overline{r}$ and $r^\star$, i.e., $ \frac{f(r^\star)-f(\overline{r}_{CLC2b})}{f(r^\star)}\times 100$ shown in Figure \ref{util}.
\end{itemize} 
We can see from the graphs that the queue occupancy under CLC2b algorithm is growing with $V$ linearly (as expected by \cite{Now,phd}) while the utility associated to the CLC2b algorithm is converging to the optimal point like $\frac{1}{V}$. The solution of CLC2b algorithm for different $V$ parameters are also shown in Figure \ref{ex-cap} (shown by the various asterisks). 
We can observe that as $V$ gets larger the solution of CLC2b algorithm converges to the optimal point.

\subsection{Stability Region for Fluid Model MQMS Systems with Stationary Continuous Channel Distribution}

We will consider a time slotted fluid model MQMS system with stationary channel distribution. 
In this case, the amount of work that arrives into (and departs from) the queues is considered to be a continuous process. We also assume that the channel state of the link from each queue to each server is modeled by a continuous random variable.
The channel state matrix in the fluid model MQMS system is defined as $C(t)=\left( C_{n,k}(t)\right),n\in \mathcal{N}, k \in \mathcal{K}, C_{n,k}(t)\in \mathbb{R}_+$.
We assume that the channel process follows a stationary distribution $f_{C(t)}(c)$ where $C(t),c \in \mathbb{R}_+^{N\times K}$; i.e., $f_C(\cdot)$ is the joint distribution of all $C_{n,k}(t)$ variables.
We can easily check that in the fluid model MQMS system, set $\widehat{V}$ is equivalent to $\mathbb{R}_+^N$ as in the fluid model MQMS system, set $\cal M$ is replaced with $\mathbb{R}_+$. Therefore, stability region for the fluid model system is characterized by the following set of inequalities.
\begin{eqnarray}
\label {cont}
\lefteqn{ \alpha {\lambda}^{\mathrm{T}}  \leq \int_{c_{n,k}=0}^{\infty}\int_{c_{n,k-1}=0}^{\infty}  \cdots \int_{c_{1,1}=0}^{\infty}    
  \max_{I \in \mathcal I}  \left(\alpha (c \circledast I )\underline{1}_K^{\mathrm{T}} \right)}
\notag ~~~~~~~~~\\  
 && \times f_{C(t)}(c)~dc_{1,1}~\cdots ~dc_{n,k-1}~dc_{n,k} ~~~ \alpha \in \mathbb{R}_+^N~~~~~~ 
\end{eqnarray}

Note that in order to characterize the stability region of the fluid model system we need to compute an infinite number of nested integrals in (\ref{cont}). In fact, the stability region of fluid model system is characterized by an infinite number of half spaces, hence the stability region is a \textit{convex surface}. In this case, depending on the channel distribution and dimension of the system, we may characterize the stability region by a finite number of non-linear inequalities instead of infinite number of linear inequalities as we show in the following example.

\el
Consider a fluid model MQMS system with two queues and one server. Assume that the channel state variables $C_{1,1}(t)$ and $C_{2,1}(t)$ are independent and follow exponential distribution with means $\mu_1$ and $\mu_2$, respectively. Such a model may be used for slow Rayleigh fading channels in low SNR regimes where SNR follows exponential distribution and the approximation $\log(1+x)\simeq x$ is used for small positive $x$. According to (\ref{cont}), the stability region is characterized by 
\begin{eqnarray}
\label {cont-ex} \alpha_1 \lambda_1 + \alpha_2 \lambda_2 \leq \frac{1}{\mu_1 \mu_2}\int_{c_{1,1}=0}^{\infty}\int_{c_{2,1}=0}^{\infty}      
 \max\{\alpha_1 c_{1,1},\alpha_2 c_{2,1}\} e^{-\frac{c_{1,1}}{\mu_1}}e^{- \frac{c_{2,1}}{\mu_2}}~dc_{2,1}dc_{1,1}~~~ \forall \alpha_1 , \alpha_2 \in \mathbb{R}_+ ~~~
\end{eqnarray}
We can write the right hand side of (\ref{cont-ex}) as
\begin{eqnarray}
\lefteqn{ \frac{1}{\mu_1 \mu_2}\int_{c_{1,1}=0}^{\infty}\int_{c_{2,1}=0}^{\frac{\alpha_1}{\alpha_2}c_{1,1}}      
\alpha_1 c_{1,1} e^{-\frac{c_{1,1}}{\mu_1}}e^{- \frac{c_{2,1}}{\mu_2}}~dc_{2,1}dc_{1,1} } \nonumber \\
&&~~~~+\frac{1}{\mu_1 \mu_2}\int_{c_{1,1}=0}^{\infty}\int_{c_{2,1}=\frac{\alpha_1}{\alpha_2}c_{1,1}}^{\infty}      
\alpha_2 c_{2,1} e^{-\frac{c_{1,1}}{\mu_1}}e^{- \frac{c_{2,1}}{\mu_2}}~dc_{2,1}dc_{1,1}  \nonumber \\
&& =
\frac{\alpha_1}{\mu_1 \mu_2}\int_{c_{1,1}=0}^{\infty} c_{1,1} e^{-\frac{c_{1,1}}{\mu_1}} \int_{c_{2,1}=0}^{\frac{\alpha_1}{\alpha_2}c_{1,1}}      
e^{- \frac{c_{2,1}}{\mu_2}}~dc_{2,1}dc_{1,1} \nonumber  \\
&&~~~~+\frac{\alpha_2}{\mu_1 \mu_2}\int_{c_{1,1}=0}^{\infty} e^{-\frac{c_{1,1}}{\mu_1}} \int_{c_{2,1}=\frac{\alpha_1}{\alpha_2}c_{1,1}}^{\infty}      
c_{2,1} e^{- \frac{c_{2,1}}{\mu_2}}~dc_{2,1}dc_{1,1}  \nonumber \\
&&=\alpha_1 \left(\mu_1- \frac{\alpha_2^2  \mu_2^2 \mu_1}{(\alpha_1 \mu_1 + \alpha_2 \mu_2)^2}  \right)
+
\alpha_2 \left( \frac{\alpha_1 \alpha_2  \mu_2^2 \mu_1}{(\alpha_1 \mu_1 + \alpha_2 \mu_2)^2} + \frac{\alpha_2 \mu_2^2}{\alpha_1 \mu_1 + \alpha_2 \mu_2}     \right)
\end{eqnarray}

Therefore, all the ordered pairs $(\lambda_1,\lambda_2)=\left(\mu_1- \frac{\alpha_2^2  \mu_2^2 \mu_1}{(\alpha_1 \mu_1 + \alpha_2 \mu_2)^2}  , \frac{\alpha_1 \alpha_2 \mu_2^2 \mu_1}{(\alpha_1 \mu_1 + \alpha_2 \mu_2)^2} + \frac{\alpha_2 \mu_2^2}{\alpha_1 \mu_1 + \alpha_2 \mu_2}  \right)$ characterize the boundary of the stability region. However, we can write $\lambda_2$ based on $\lambda_1$, $\mu_1$ and $\mu_2$ as follows.
\begin{eqnarray}
\mu_1- \frac{\alpha_2^2  \mu_2^2 \mu_1}{(\alpha_1 \mu_1 + \alpha_2 \mu_2)^2} =\lambda_1 \nonumber
\end{eqnarray}
So, 
\begin{eqnarray}
\frac{\alpha_2  \mu_2 }{\alpha_1 \mu_1 + \alpha_2 \mu_2}=\sqrt{1-\frac{\lambda_1}{\mu_1}} 
~~~~~~~~\text{and}~~~~~~~
\frac{\alpha_1  \mu_1}{\alpha_1 \mu_1 + \alpha_2 \mu_2}=1-\sqrt{1-\frac{\lambda_1}{\mu_1}} \nonumber
\end{eqnarray}
Therefore,
\begin{eqnarray}
\label{stab-2q}\lambda_2 =\mu_2 \left( \sqrt{1-\frac{\lambda_1}{\mu_1}}\right) \left( 2-\sqrt{1-\frac{\lambda_1}{\mu_1}}  \right)
\end{eqnarray}
\begin{figure}[tp]
    \centering
    \includegraphics[width=0.7\textwidth]{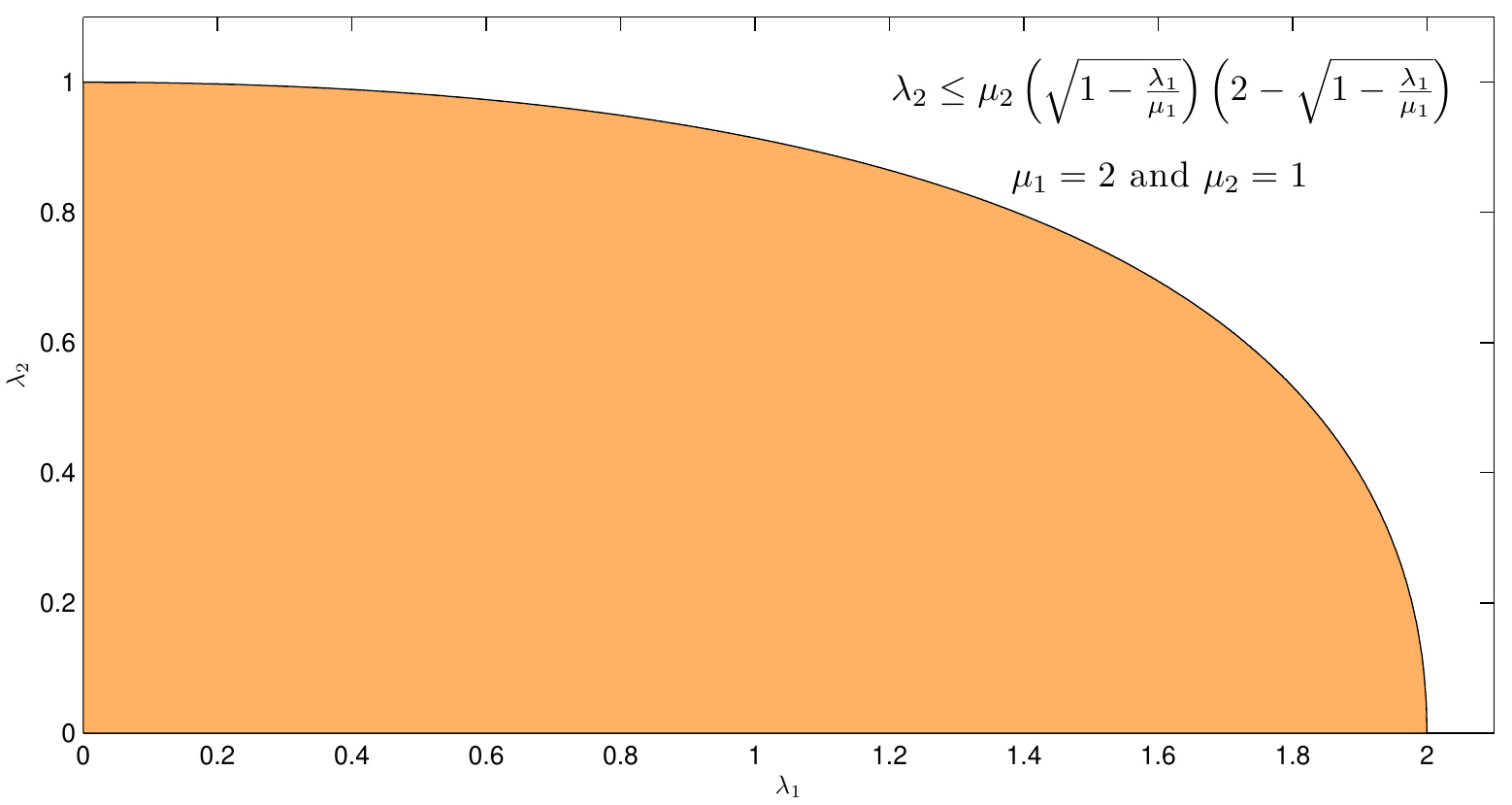}
	\caption{ Stability region for $\mu_1=2$ and $\mu_2 =1$}
	\label{cont-2q} \vspace{-2mm}
\end{figure}

As we can see from equation (\ref{stab-2q}), the stability region in this example is characterized by just one non-linear inequality $\lambda_2 \leq \mu_2 \left( \sqrt{1-\frac{\lambda_1}{\mu_1}}\right) \left( 2-\sqrt{1-\frac{\lambda_1}{\mu_1}}  \right)$ and two linear inequalities $\lambda_1 \geq 0$ and $\lambda_2 \geq 0$.
The stability region for this example for $\mu_1=2$ and $\mu_2 =1$ is illustrated in Figure \ref{cont-2q}. The characterization of the general fluid model stability region is beyond the scope of the paper and left as a possible future research problem.

\section{Conclusions and Future Work}\label{conc}
In this paper, we introduced a linear algebraic representation of the network stability region (capacity region) polytope of multi-queue multi-server (MQMS) queueing system with stationary channel distribution and stationary arrival processes. To this end, we obtained the necessary and sufficient conditions for the stability of the system for a general arrival process with finite first and second moments. For stationary arrival processes, we showed that these conditions establish the network stability region of the system given by (\ref{t1_2}). For the stability region polytope, we explicitly determined all the coefficients $\alpha \in \widehat{V}$ of all the half spaces which are required to characterize the stability region polytope. We also argued that in general it may be (computationally) hard to quantify set $\widehat{V}$. In this case, although we may add some redundant inequalities, we can use a superset of $\widehat{V}$, namely the set $\mathcal{W}^N-\{ \underline{0}_N\}$ ($\widehat{V} \subseteq \mathcal{W}^N- \{ \underline{0}_N\}$) instead of $\widehat{V}$ in (\ref{t1_2}).
An upper bound was also obtained for the average queueing delay of Maximum Weight (MW) server allocation policy which is a throughput optimal policy for MQMS system. We finally considered the stability region for a fluid model MQMS. In this case, we determine the stability region by an infinite set of linear inequalities given by (\ref{cont}). By use of an example we showed that depending on the channel distribution and the number of queues, we may characterize the stability region by a finite set of non-linear inequalities instead of infinite number of linear inequalities. However, the general problem of stability region characterization for the fluid model MQMS can be considered as a possible future research problem.

\appendices
\section{Proof of Lemma \ref{l6}} \label{a1}
\begin{IEEEproof}
As we explained before, ${\cal I}^{\alpha}$ denotes a set of allocation matrices of $\{I^{\alpha}_s, s\in \cal S\}$ that maximize the right hand side of (\ref{l4_1}) and we discussed that ${\cal I}^{\alpha}$ is not a unique solution of (\ref{argmax}) and there may be more than one set of allocation matrices of ${\cal I}^{\alpha}$ whose elements maximize (\ref{argmax}). Let $\mathbb{I}^{\alpha}=\{{\cal I}^{\alpha}_{i}, ~1\leq i \leq |\mathbb{I}^{\alpha}|\}$ denote the set of all distinguished solutions of (\ref{argmax}). Obviously $|\mathbb{I}^{\alpha}|<\infty$, since set $\cal I$ is finite. Note that each solution ${\cal I}^{\alpha}_{i}=\{ I_{s,i}^{\alpha} , s\in {\cal S}  \}$, $1\leq i \leq |\mathbb{I}^{\alpha}|$ is corresponding to a deterministic policy with the average transmission rate vector  
\begin{eqnarray}
\label{l6_1}
\label{r-alpha}R^{\alpha}_i={\left( \displaystyle\sum_{s\in {\cal S}}\pi_s \left(C_s\circledast I_{s,i}^{\alpha} \right)\underline{1}_K^{\mathrm{T}}\right)}^{\mathrm{T}}
\end{eqnarray}
and therefore, according to (\ref{conv}) each solution of (\ref{argmax}) is associated with a vertex of polytope $\cal P$. Since there are $|\mathbb{I}^{\alpha}|$ distinguished solution of (\ref{argmax}), $|\mathbb{I}^{\alpha}|$ vertices of polytope $\cal P$ are on the hyperplane associated to half space (\ref{l4_1}), i.e.,
\begin{eqnarray} \label{l4_1_hyp}
\alpha \overline{E[A(t)]^{\mathrm{T}}} = \displaystyle\sum_{s\in \mathcal{S}}\pi_s 
  \max_{I \in \mathcal I}  \left(\alpha (C_s \circledast I )\underline{1}_K^{\mathrm{T}} \right)
\end{eqnarray}
Therefore, each face of $\cal P$ defined by hyperplane (\ref{l4_1_hyp}) is represented by $F_{\alpha}=conv.hull_{1\leq i \leq |\mathbb{I}^{\alpha}|} R^{\alpha}_i$.
Since polytope $\cal P$ is full dimensional, according to Definition \ref{facet}, each facet of $\cal P$ must be of dimension $N-1$. In the following, we will show that the dimension of face $F_{\alpha}$ is less than $N-1$ if $\alpha \in {\mathbb R^N_+}-V$ and therefore the hyperplane (\ref{l4_1_hyp}) associated to such $\alpha$ is not a facet defining hyperplane.

Consider a vector $\alpha\in  {\mathbb R^N_+}-V$. For $\alpha$ we have the following property:
\begin{eqnarray}
\label{l6_3}
\lefteqn{\exists ~(\mathcal{U} \subset{\cal N}, \mathcal{U}\neq\varnothing, \alpha_{\mathcal{U}}\neq \underline{0}_{|\mathcal{U}|}, \alpha_{\mathcal{U}^c}\neq \underline{0}_{|\mathcal{U}^c|}) :} \nonumber \\ 
&& \forall ~ (i  \in  \mathcal{U} , j \in  \mathcal{U}^c ,~m,n \in \mathcal M , \alpha_{i}, \alpha_{j}, m,n\neq 0 )~ ~\alpha_{i}m \neq \alpha_{j}n
\end{eqnarray}
In other words, there exists a partitioning of vector $\alpha$ into two (sub)vectors $\alpha_{{\cal U}}$ and $\alpha_{{\cal U}^c}$ such that no non-zero element of $\alpha_{{\cal U}}$ is proportional to no non-zero element of $\alpha_{{\cal U}^c}$ by the ratio of two non-zero elements in $\cal M$. Assume that ${\cal U}_0$ denotes the subset $\cal U$ satisfying (\ref{l6_3}).
Now, consider inequality (\ref{l4_1}). The right hand side of (\ref{l4_1}) can also be expressed as 
\begin{eqnarray}
\label{l6_2}
\displaystyle\sum_{s\in \mathcal{S}}\pi_s 
  \max_{I \in \mathcal I}  \left(\alpha (C_s \circledast I )\underline{1}_K^{\mathrm{T}} \right)=\displaystyle\sum_{s\in \mathcal{S}}\pi_s 
 \displaystyle\sum_{k=1}^K \max_{J \in \mathcal J}  \left\langle \alpha , ( (C_s^{\downarrow k })^{\mathrm{T}} \circledast J )\right\rangle,
\end{eqnarray}
where $C_s^{\downarrow k }$ is the $k$'th column of matrix $C_s$ and $\cal J$ is the set of all binary vectors of size $N$ with $\left\langle J, \underline{1}_N\right\rangle=1$, i.e.,
\begin{eqnarray}
{\cal J}:=\{J=(J_1,J_2,...,J_N)\mid J_n \in \{0,1\}~\forall n \in {\cal N},\left\langle J, \underline{1}_N\right\rangle=1\}.
\end{eqnarray}
Each $J_s^k(\alpha)=\arg\max_{J \in \mathcal J}  \left\langle\alpha ,((C_s^{{\downarrow k}})^{\mathrm{T}} \circledast J )\right\rangle$ is corresponding to the $k$'th column of an $I_s^{\alpha}$ in some ${\cal I}^{\alpha} \in \mathbb{I}^{\alpha}$. Note that $\max_{J \in \mathcal J} \left\langle\alpha ,((C_s^{{\downarrow k}})^{\mathrm{T}} \circledast J )\right\rangle \geq 0$. There are $K|{\cal S}|$
of such maximization terms in the right hand side of (\ref{l6_2}). 
 Each of the $K|{\cal S}|$ maximization terms may have multiple solutions (non-unique solutions). Suppose that $\psi_s^k(\alpha)$ denotes the set of all distinct solutions to $J_s^k(\alpha)=\arg\max_{J \in \mathcal J}  \left\langle\alpha ,((C_s^{{\downarrow k}})^{\mathrm{T}} \circledast J )\right\rangle$, i.e.,
\begin{eqnarray}
\psi_s^k(\alpha)=\{J_s^k(\alpha)=\arg\max_{J \in \mathcal J}  \left\langle\alpha ,((C_s^{{\downarrow k}})^{\mathrm{T}} \circledast J )\right\rangle~\forall s\in \mathcal{S},~k \in \mathcal{K}\}
\end{eqnarray}

For a given $ s\in \mathcal{S},~k \in \mathcal{K} $ with $\left\langle\alpha ,\left((C_s^{{\downarrow k}})^{\mathrm{T}} \circledast {J_s^k(\alpha)} \right)\right\rangle > 0$, we can easily observe that for all the elements of $ \psi_s^k(\alpha) $ either we have $\sum_{u \in {\cal U}_0} J_{s,u}^k(\alpha)=1$ or $\sum_{u \in {\cal U}_0^c}J_{s,u}^k(\alpha)=1$. This result comes directly from (\ref{l6_3}) and states that for a given $ s\in \mathcal{S},~k \in \mathcal{K} $, \textit{for all} the solutions of $J_s^k(\alpha)=\arg\max_{J \in \mathcal J}  \left\langle\alpha ,((C_s^{{\downarrow k}})^{\mathrm{T}} \circledast J )\right\rangle$ with $\max_{J \in \mathcal J}  \left\langle\alpha ,((C_s^{{\downarrow k}})^{\mathrm{T}} \circledast J )\right\rangle>0$, the index associated to element ``1" of $J_s^k(\alpha)$ is either in ${\cal U}_0$ or ${\cal U}_0^c$.  
In other words, it is not possible to have the index of element ``1" in ${\cal U}_0$ for some solutions and in ${\cal U}_0^c$ for other solutions of $J_s^k(\alpha)=\arg\max_{J \in \mathcal J}  \left\langle\alpha ,((C_s^{{\downarrow k}})^{\mathrm{T}} \circledast J )\right\rangle$. 
Using this property, we can partition set of all $\psi_s^k(\alpha)$ with $\left\langle\alpha ,((C_s^{{\downarrow k}})^{\mathrm{T}} \circledast {J_s^k(\alpha)} )\right\rangle > 0$ into two disjoint subsets. We denote the set of all solution sets $\psi_s^k(\alpha)$ with $\left\langle\alpha ,((C_s^{{\downarrow k}})^{\mathrm{T}} \circledast J_s^k (\alpha) )\right\rangle > 0$ and $\sum_{u \in {\cal U}_0} J_{s,u}^k(\alpha)=1$ as ${\cal A}_\alpha$ and set of all solution sets $\psi_s^k(\alpha)$ with $\left\langle\alpha ,((C_s^{{\downarrow k}})^{\mathrm{T}} \circledast J_s^k(\alpha) )\right\rangle > 0$ and $\sum_{u \in {\cal U}_0^c} J_{s,u}^k(\alpha)=1$ as ${\cal B}_\alpha$.
We introduce $N$ dimensional vectors $\alpha'_{{\cal U}_0}$ and $\alpha'_{{\cal U}^c_0}$ as follows.
\[\alpha'_{{\cal U}_0}=(\alpha'_1,\alpha'_2,...,\alpha'_N):~\alpha'_n= \left\{
  \begin{array}{l l}
    \alpha_n & \quad \text{if $n \in \mathcal{U}$ }\\
    0 & \quad \text{if $n \notin \mathcal{U}$}\\
  \end{array} \right. \]
\[\alpha'_{{\cal U}_0^c}=(\alpha'_1,\alpha'_2,...,\alpha'_N):~\alpha'_n= \left\{
  \begin{array}{l l}
    \alpha_n & \quad \text{if $n \in \mathcal{U}^c$ }\\
    0 & \quad \text{if $n \notin \mathcal{U}^c$}\\
  \end{array} \right. \]
 Now, let us introduce the following two hyperplanes.
\begin{eqnarray}
\label{l6_4}\alpha'_{{\cal U}_0}\overline{E[{A_{\cal U}(t)}]^{\mathrm{T}}}= \displaystyle\sum_{k,s: \psi_s^k(\alpha) \in {\cal A}_\alpha} \pi_s \max_{J \in \mathcal{J}} \left\langle\alpha ,((C_s^{{\downarrow k}})^{\mathrm{T}} \circledast J )\right\rangle ~
 \\
\label{l6_5}\alpha'_{{\cal U}_0^c}\overline{E[A_{{\cal U}^c}(t)]^{\mathrm{T}}}= \displaystyle\sum_{k,s: \psi_s^k(\alpha) \in {\cal B}_\alpha} \pi_s \max_{J \in \mathcal{J}} \left\langle\alpha ,((C_s^{{\downarrow k}})^{\mathrm{T}} \circledast J )\right\rangle
\end{eqnarray}

In the following, we will prove that all the vertices $R^{\alpha}_i, ~1\leq i \leq |\mathbb{I}^{\alpha}|$ in (\ref{r-alpha}) satisfy both (\ref{l6_4}) and (\ref{l6_5}). In other words, all the vertices $R^{\alpha}_i$ are located on both hyperplanes (\ref{l6_4}) and (\ref{l6_5}) and therefore on the intersection of both. The hyperplanes of (\ref{l6_4}) and (\ref{l6_5}) are in an $N$ dimensional space whose dimensions are at most $N-1$. Thus, dimension of their intersection is at most $N-2$. In other words, $\dim F_{\alpha}= \dim ( conv.hull_{1\leq i \leq |\mathbb{I}^{\alpha}|} R^{\alpha}_i)<N-1$. Therefore $F_\alpha$ is not a facet for $\cal P$.

Note that columns of any allocation matrix $I^{\alpha}_s$ are vectors $J_s^k(\alpha) \in \psi_s^k(\alpha)$ for all $k \in \cal K$.
Consider a particular transmission rate vector $ R^{\alpha}_i $ (equation (\ref{r-alpha})) located on hyperplane (\ref{l4_1_hyp}). As we discussed before, each allocation matrix $I^{\alpha}_{s,i}$ in (\ref{r-alpha}) is created by concatenation of a set of vectors $J_{s,i}^k(\alpha) \in \psi_s^k(\alpha)$ as the columns of $I^{\alpha}_{s,i}$. Now, we can check that all the vertices $R^{\alpha}_i, ~1\leq i \leq |\mathbb{I}^{\alpha}|$ in (\ref{r-alpha}) satisfy both (\ref{l6_4}) and (\ref{l6_5}) as follows. 
\begin{eqnarray}
\lefteqn{ \alpha'_{{\cal U}_0}{\left( \displaystyle\sum_{s\in {\cal S}}\pi_s \left(C_s\circledast I_{s,i}^{\alpha} \right)\underline{1}_K^{\mathrm{T}}\right)}} \nonumber\\
&&= \displaystyle\sum_{s\in \mathcal{S}}\pi_s 
 \displaystyle\sum_{k=1}^K  \left\langle \alpha'_{{\cal U}_0} , ( (C_s^{\downarrow k })^{\mathrm{T}} \circledast J_{s,i}^k (\alpha) )\right\rangle =\displaystyle\sum_{k,s: \psi_s^k(\alpha) \in {\cal A}_\alpha} \pi_s \max_{J \in \mathcal{J}} \left\langle\alpha ,((C_s^{{\downarrow k}})^{\mathrm{T}} \circledast J )\right\rangle
\end{eqnarray}
and
\begin{eqnarray}
\lefteqn{ \alpha'_{{\cal U}_0^c} {\left( \displaystyle\sum_{s\in {\cal S}}\pi_s \left(C_s\circledast I_{s,i}^{\alpha} \right)\underline{1}_K^{\mathrm{T}}\right)}} \nonumber\\
&&= \displaystyle\sum_{s\in \mathcal{S}}\pi_s 
 \displaystyle\sum_{k=1}^K  \left\langle \alpha'_{{\cal U}_0^c} , ( (C_s^{\downarrow k })^{\mathrm{T}} \circledast J_{s,i}^k(\alpha) )\right\rangle =\displaystyle\sum_{k,s: \psi_s^k(\alpha) \in {\cal B}_\alpha} \pi_s \max_{J \in \mathcal{J}} \left\langle\alpha ,((C_s^{{\downarrow k}})^{\mathrm{T}} \circledast J )\right\rangle
\end{eqnarray}
 Thus, all the vertex (\ref{r-alpha}) satisfy both (\ref{l6_4}) and (\ref{l6_5}) and therefore their intersection. This shows that $F_\alpha$ has dimension less than $N-1$ and therefore is not a facet of $\cal P$. 
\end{IEEEproof}

\section{Proof of Lemma \ref{l5}}\label{a2}
\begin{IEEEproof}
 Consider vector $\alpha=(\alpha_1,\alpha_2,...,\alpha_N) \in V$. From (\ref{Vdef}), it is obvious that any vector derived by the multiplication of a positive scalar to $\alpha$ also belongs to $V$, i.e.,
\begin{eqnarray}\label{l5_2}
\alpha \in V \Longrightarrow q\alpha \in V~~ \forall q\in \mathbb{R}_+
\end{eqnarray}

Also note that vectors $\alpha$ and $q\alpha$ will result in the same inequalities in (\ref{l4_1}) and therefore the same face defining hyperplanes of polytope $\cal P$. Thus, we say that $\alpha$ and $q\alpha$ are \textit{equivalent} and we write $\alpha\equiv q\alpha$.

Consider a vector $\alpha \in V$. Vector $\alpha$ may be a zero vector in which case it does not contributes to  a face of $\cal P$ as the inequality (\ref{l4_1}) will result in an obvious equality $0=0$. Hence, we assume that $\alpha$ is not a zero vector. 

Now, we perform the following process on vector $\alpha$.
Pick a non-zero elements of $\alpha$, let say $\alpha_{e_1}$ and form $\alpha_{\{e_1\}}$ and  $\alpha_{{\cal N}- \{e_1\}}$. The size of these vectors are 1 and $N-1$, respectively. According to (\ref{Vdef}), the vector $\alpha_{{\cal N}- \{e_1\}}$ should be either a zero vector in which case we stop or there should be a non-zero element $\alpha_{e_2}$ in $\alpha_{{\cal N}- \{e_1\}}$ and non-zero $m_1,n_1\in \cal M$ such that $\alpha_{e_2}=\alpha_{e_1} \frac{m_1}{n_1}$. If so, we will continue and form $\alpha_{\{e_1,e_2\}}$ and $\alpha_{{\cal N}- \{e_1,e_2\}}$. According to (\ref{Vdef}), the vector $\alpha_{{\cal N}- \{e_1,e_2\}}$ is either a zero vector in which case we will stop or there should exist a non-zero $\alpha_{e_3}\in \alpha_{{\cal N}- \{e_1,e_2\}}$ and non-zero $m_2,n_2$ such that either $\alpha_{e_3}=\alpha_{e_1} \frac{m_2}{n_2}$ or $\alpha_{e_3}=\alpha_{e_2} \frac{m_2}{n_2}=\alpha_{e_1} \frac{m_1}{n_1} \frac{m_2}{n_2}$. We repeat this procedure until either for some $\nu<N$ the vector $\alpha_{{\cal N}-\{e_1,e_2,...,e_{\nu}\}}$ is a zero vector or we form vector $\alpha_{\{e_1,e_2,...,e_N\}}$ in which case $\nu=N$. In both cases, we can see that all the non-zero elements of $\alpha$ may be expressed by the multiplication of $\alpha_1$ and the rational numbers derived by non-zero elements of set $\cal M$. By dividing the vector $\alpha$ by scalar $\alpha_{e_1}$ and multiplying it by $\prod_{j=1}^{\nu-1}n_j$, we will obtain an equivalent vector for $\alpha$ named $\beta$ whose elements belong to the set of all integers which were derived by multiplication of $N-1$ elements of set $\cal M$, i.e.,
\begin{eqnarray}
\label{l5_3}
\alpha \equiv \beta~,~\beta_n \in {\cal W}=\left\{z \in \mathbb{Z_+} \mid z=\prod_{j=1}^{N-1}m_j,~m_j \in \cal M \right\}~~\forall n \in \{1,2,...,N\}
\end{eqnarray}  
 
Therefore, each vector $\alpha \in V$ is equivalent to a vector $\beta$ whose elements come from set $\cal W$.
We define the set of all $\beta$ vectors with the above property as $\widehat{V}$. Note that $\widehat{V} \subseteq \mathcal{W}^N$ and is a finite set since $\cal W$ is finite.
Thus, for each vector $\alpha \in V$ there exists an equivalent vector $\beta \in \widehat{V}$ and the result follows.
\end{IEEEproof}

\section{Proof of Theorem \ref{sc}}\label{a3}
\begin{IEEEproof} 
We will start with the Lyapunov function evaluation. we will use the quadratic function (\ref{lypfunc}) as our Lyapunov function. The Lyapunov drift for two successive time slots has the following form.
\begin{eqnarray}
\label{sc3}
\lefteqn {D(t+1)=E[L(X(t+1))-L(X(t)) \mid X(t)]} \nonumber \\ & &= E\left[\displaystyle\sum_{n=1}^N X_n^2(t+1)-X_n^2(t) \mid X(t)\right] =  E\left[\displaystyle\sum_{n=1}^N (X_n(t+1)-X_n(t))^2 \mid X(t)\right]\nonumber \\ \label{sc3}
& &+ {\;\;\;}2 E\left[\displaystyle\sum_{n=1}^N X_n(t)(X_n(t+1)-X_n(t))\mid X(t)\right]
\end{eqnarray}

For the the first term we have:
\begin{eqnarray}
\label{sc4}
\notag \lefteqn{E\left[\displaystyle\sum_{n=1}^N (X_n(t+1)-X_n(t))^2 \mid X(t)\right]}\\ \notag
& &= E\left[\displaystyle\sum_{n=1}^N (A_n(t+1)-\displaystyle\sum_{k=1}^K H_{n,k}(t+1))^2 \mid X(t)\right]= E\left[\displaystyle\sum_{n=1}^N A_n^2(t+1) \mid X(t)\right] \\ \nonumber & & -2E\left[\displaystyle\sum_{n=1}^N\displaystyle\sum_{k=1}^K A_n(t+1)H_{n,k}(t+1) \mid X(t)\right] + \;\;\; E\left[\displaystyle\sum _{n=1}^N \left ( \displaystyle\sum_{k=1}^K H_{n,k}(t+1) \right )^2 \mid X(t)\right] \nonumber \\
&& \leq N A_{max}^2 + \displaystyle\sum _{n=1}^N E\left[ \left ( \displaystyle\sum_{k=1}^K H_{n,k}(t+1) \right )^2 \mid X(t)\right]
\end{eqnarray}

For the second term in (\ref{sc3}) we have
\begin{eqnarray}
\label{sc7}\notag \lefteqn {E\left[\displaystyle\sum_{n=1}^N X_n(t)(X_n(t+1)-X_n(t))\mid X(t)\right]}\\ 
\notag &&= E\left[\displaystyle\sum_{n=1}^N X_n(t) \left( A_n(t+1)-\displaystyle\sum_{k=1}^K H_{n,k}(t+1) \right)\mid X(t)\right] \\
&& = \displaystyle\sum_{n=1}^N X_n(t) E [A_n(t+1)]   -\displaystyle\sum_{n=1}^N E \left[ X_n(t) \displaystyle\sum_{k=1}^K H_{n,k}(t+1) \mid X(t)\right] 
\end{eqnarray}
Therefore, the Lyapunov drift $D(t+1)$ can be bounded by 
\begin{eqnarray} \label{sc7.5}
\lefteqn{D(t+1)\leq N A_{max}^2 + \displaystyle\sum _{n=1}^N E\left[ \left ( \displaystyle\sum_{k=1}^K H_{n,k}(t+1) \right )^2 \mid X(t)\right]} \nonumber \\
&&~~~~~~~~~~ +2\displaystyle\sum_{n=1}^N X_n(t) E [A_n(t+1)]   -2\displaystyle\sum_{n=1}^N E \left[ X_n(t) \displaystyle\sum_{k=1}^K H_{n,k}(t+1) \mid X(t)\right] 
\end{eqnarray}
In the following, we show that the Lyapunov drift in (\ref{sc7.5}) is bounded as follows.
\begin{eqnarray} \label{sc7.75}
\lefteqn{D(t+1)\leq N A_{max}^2 + \displaystyle\sum _{n=1}^N E\left[ \left ( \displaystyle\sum_{k=1}^K C_{n,k}(t+1)I_{n,k}(t+1) \right )^2 \mid X(t)\right]} \nonumber \\
&&~~~~~ +2\displaystyle\sum_{n=1}^N X_n(t) E [A_n(t+1)]   -2\displaystyle\sum_{n=1}^N E \left[ X_n(t) \displaystyle\sum_{k=1}^K C_{n,k}(t+1)I_{n,k}(t+1) \mid X(t)\right] 
\end{eqnarray}
To prove (\ref{sc7.75}), note that for each queue $n\in \cal N$ one of the following conditions is satisfied.
\begin{itemize}
\item $\sum_{k=1}^K H_{n,k}(t+1) =\sum_{k=1}^K   C_{n,k}(t+1)I_{n,k}(t+1)$: In this case, we can easily check that for all $n \in \cal N$
\begin{eqnarray} \label{sc7.8}
\lefteqn{E[X_n^2(t+1)-X_n^2(t)\mid X(t)]\leq N A_{max}^2 } \nonumber \\
&& + E\left[ \left( \displaystyle\sum_{k=1}^K C_{n,k}(t+1)I_{n,k}(t+1) \right)^2 \mid X(t)\right] 
 \nonumber \\
&& +2 X_n(t) E [A_n(t+1)]- 2E \left[ X_n(t) \displaystyle\sum_{k=1}^K C_{n,k}(t+1)I_{n,k}(t+1) \mid X(t)\right] 
\end{eqnarray}

\item $\sum_{k=1}^K H_{n,k}(t+1) < \sum_{k=1}^K C_{n,k}(t+1)I_{n,k}(t+1)$: In this case, $\sum_{k=1}^K H_{n,k}(t+1) = X_n(t)$ and $ X_n(t+1)=A_n(t+1) $ as there are not enough packets in queue $n$ to be server at time slot $t+1$. It is not hard to check that inequality (\ref{sc7.8}) is also satisfied in this case.
\end{itemize} 
According to the above discussion, we can observe that inequality (\ref{sc7.8}) is satisfied for all the queues and therefore (\ref{sc7.75}) follows. 

Using the the fact that \(\displaystyle\sum_{k=1}^K C_{n,k}(t+1)I_{n,k}(t+1) \geq 0\), we get the following inequality.
\begin{eqnarray}
\label{sc5} & \displaystyle\sum _{n=1}^N \left ( \displaystyle\sum_{k=1}^K C_{n,k}(t+1)I_{n,k}(t+1) \right )^2  \leq \left ( \displaystyle\sum_{n=1}^N\displaystyle\sum_{k=1}^K C_{n,k}(t+1)I_{n,k}(t+1)\right)^2 \leq (MK)^2 
\end{eqnarray}
Hence, the Lyapunov drift (\ref{sc7.75}) is bounded by 
\begin{eqnarray} \label{sc7.9}
\lefteqn{D(t+1)\leq N A_{max}^2 + (MK)^2+2\displaystyle\sum_{n=1}^N X_n(t) E [A_n(t+1)]  }\nonumber \\ &&~~~~~~~~~~~~-2\displaystyle\sum_{n=1}^N E \left[ X_n(t) \displaystyle\sum_{k=1}^K C_{n,k}(t+1)I_{n,k}(t+1) \mid X(t)\right] 
\end{eqnarray} 
By conditioning the last term of (\ref{sc7.9}) on the channel state at time slot $t+1$, we will have
\begin{eqnarray}
\label{sc8}\notag \lefteqn{D(t+1)\leq N A_{max}^2 + (MK)^2} \\
 && +2 E \left[ X(t) \left(  A^{\mathrm{T}}(t+1)-\sum_{s \in \mathcal{S} } \pi_s  \left( C_s \circledast I(t+1)\right)\underline{1}_K^{\mathrm{T}} \right)  \mid  X(t) \right]
\end{eqnarray}

Note that the allocation matrix $I$ in (\ref{sc8}) depends on the selected policy. According to (\ref{mw}), we can see that by selecting MW policy, the second term of (\ref{sc8}) will be minimized and therefore the right hand side term in (\ref{sc8}) will be minimized over all the existing server allocation policies.
In other words, for MW policy and any arbitrary server allocation policy $\Delta$ we have
\begin{eqnarray}
\label{sc9}\notag \lefteqn{D^{MW}(t+1)\leq N A_{max}^2 + (MK)^2} \\
 && ~~+2 E \left[ X(t) \left(  A^{\mathrm{T}}(t+1)-\sum_{s \in \mathcal{S} } \pi_s  \left( C_s \circledast I^{MW}(t+1)\right)\underline{1}_K^{\mathrm{T}} \right)  \mid  X(t) \right]
 \nonumber \\
 && \leq N A_{max}^2 + (MK)^2 \nonumber\\
 && ~~+2 E \left[ X(t) \left(  A^{\mathrm{T}}(t+1)-\sum_{s \in \mathcal{S} } \pi_s  \left( C_s \circledast I^{\Delta}(t+1)\right)\underline{1}_K^{\mathrm{T}} \right)  \mid  X(t) \right]
\end{eqnarray}
where $I^{MW}(t+1)$ and $I^{\Delta}(t+1)$ are the allocation matrices of policies MW and $\Delta$, respectively and $D^{MW}(t+1)$ is the Lyapunov drift of policy MW at time slot $t+1$.

It is important to note that the set of inequalities in (\ref{sc1}) determines an open polytope $\mathcal{P}'$ for which we have $\mathcal{P}'=\mathcal{P}-bound(\mathcal{P})$. In fact, Theorem \ref{sc} states that if the vector $E[A(t)]$ at each time slot $t$ is strictly inside polytope $\mathcal{P}$ then the system is stable. Now, consider an MQMS system with arrival processes for which we have $E[A(t)] \in \mathcal{P}'$ at each time slot $t$. 
Suppose that $\delta>0$ is a positive real number such that $E[A(t)]+\delta \underline{1}_N \in bound (\mathcal{P})$. Therefore, $E[A(t)]+\delta \in \mathcal{P}~ \forall t$. According to our discussion about randomized policies, there exists a randomized policy $I^{rnd}=\{I^{rnd}(t)\}_{t=1}^{\infty}$ for which we will have
\begin{eqnarray}
\label{sc9}
E\left[ \left( \sum_{s \in \mathcal{S} } \pi_s  \left( C_s \circledast I^{rnd}(t+1)\right)\underline{1}_K^{\mathrm{T}}\right)^{\mathrm{T}}\right] =  E[A(t+1)]+ \delta \underline{1}_N 
\end{eqnarray}
Therefore, 
\begin{eqnarray} \label{sc9.5}
E \left[ X(t) \left(  A^{\mathrm{T}}(t+1)-\sum_{s \in \mathcal{S} } \pi_s  \left( C_s \circledast I^{\Delta}(t+1)\right)\underline{1}_K^{\mathrm{T}} \right)  \mid  X(t) \right]=-\delta \sum_{n=1}^N X_n(t)
\end{eqnarray}
Putting together from (\ref{sc8}) to (\ref{sc9.5}), we conclude that
\begin{eqnarray}
D^{MW}(t+1)\leq N A_{max}^2+(MK)^2 -2 \delta \sum_{n=1}^N X_n(t)
\end{eqnarray}
and according to (\ref{lyapunov}) the stability of MQMS system is demonstrated.
\end{IEEEproof}

\bibliographystyle{ieeetran}
\bibliography{Ref}

\end{document}